\documentclass[a4paper]{amsproc}
\usepackage{graphicx}

\newtheorem{theorem}{Theorem}[section]
\newtheorem{lemma}[theorem]{Lemma}
\newtheorem{corollary}[theorem]{Corollary}
\newtheorem{proposition}[theorem]{Proposition}
\newtheorem{remark}[theorem]{Remark}
\newtheorem{definition}[theorem]{Definition}

\theoremstyle{remark}

\numberwithin{equation}{section}

\newcommand{\bz}{{\mathbb B}}
\newcommand{\cz}{{\mathbb C}}
\newcommand{\dz}{{\mathbb D}}

\newcommand{\nz}{{\mathbb N}}
\newcommand{\rz}{{\mathbb R}}

\newcommand{\calA}{\mathcal{A}}
\newcommand{\calB}{\mathcal{B}}
\newcommand{\calC}{\mathcal{C}}
\newcommand{\calD}{\mathcal{D}}
\newcommand{\calE}{\mathcal{E}}
\newcommand{\calF}{\mathcal{F}}
\newcommand{\calH}{\mathcal{H}}
\newcommand{\calK}{\mathcal{K}}
\newcommand{\calL}{\mathcal{L}}
\newcommand{\calM}{\mathcal{M}}

\newcommand{\calR}{\mathcal{R}}
\newcommand{\calS}{\mathcal{S}}

\newcommand{\ulA}{\underline{A}}

\newcommand{\forget}[1]{}

\newcommand{\ci}{\mathcal{C}^\infty}
\newcommand{\cicomp}{\mathcal{C}^\infty_{\text{\rm comp}}}
\newcommand{\cig}{\mathcal{C}^{\infty,\gamma}}

\newcommand{\cii}{\mathcal{C}^{\infty,\infty}}
\newcommand{\cl}{{\rm cl}}
\newcommand{\dbar}{d\hspace*{-0.08em}\bar{}\hspace*{0.1em}}
\newcommand{\eps}{\varepsilon}

\newcommand{\intb}{\text{\rm int}\,{\mathbb B}}
\newcommand{\intd}{\text{\rm int}\,{\mathbb D}}
\newcommand{\op}{\text{\rm op}}
\newcommand{\opm}[1]{\text{\rm op}_M^{#1}}
\newcommand{\pit}{\,{\widehat{\otimes}}_\pi\,}
\newcommand{\re}{\text{\rm Re}\,}
\newcommand{\rpbar}{\overline{{\mathbb R}}_+}
\newcommand{\schnitt}{\mathop{\mbox{\Large$\cap$}}}

\newcommand{\smsum}{\mathop{\mbox{\large$\sum$}}}
\newcommand{\spalte}[2]{\mbox{\LARGE$\binom{#1}{#2}$}}
\newcommand{\spk}[1]{\left<#1\right>}
\newcommand{\st}{\mbox{\boldmath$\;|\;$\unboldmath}}
\newcommand{\trinorm}[1]%
    {|\hspace*{-1pt}|\hspace*{-1pt}|#1|\hspace*{-1pt}|\hspace*{-1pt}|}

\newcommand{\wh}{\widehat}
\newcommand{\wt}{\widetilde}
\renewcommand{\Re}{{\rm Re}\,}

\begin{document}

\title[$H_\infty$-Calculus for Cone Differential Operators]
      {Bounded $\mathbf{H_\infty}$-Calculus for Differential Operators \\
       on Conic Manifolds with Boundary}
\author{S.\ Coriasco}
\address{Universita di Torino, Dipartimento di Matematica, Via Carlo Alberto 10, 
         10123 Torino, Italy}
\author{E.\ Schrohe}
\address{Universit\"at Hannover, Institut f\"ur Mathematik,
         Welfengarten 1, 30167 Hannover, Germany}
\author{J.\ Seiler}
\address{Universit\"at Hannover, Institut f\"ur Angewandte Mathematik,
         Welfengarten 1, 30167 Hannover, Germany}
\email{sandro.coriasco@unito.it, schrohe@math.uni-hannover.de, seiler@ifam.uni-hannover.de}


\begin{abstract}
We derive conditions that ensure the existence of a bounded $H_\infty$-calculus 
in weighted $L_p$-Sobolev spaces for 
closed extensions $\underline{A}_T$ of a differential operator $A$ on a conic 
manifold with boundary, subject to differential boundary conditions $T$. In general, 
these conditions ask for a particular pseudodifferential structure of the resolvent 
$(\lambda-\underline{A}_T)^{-1}$ in a sector $\Lambda\subset\cz$. In case of the 
minimal extension they reduce to parameter-ellipticity of the boundary 
value problem $\binom{A}{T}$. Examples concern the Dirichlet and Neumann Laplacians.  
\end{abstract}

\maketitle

\tableofcontents
\section{Introduction}\label{sectionh1}

Establishing the existence of a bounded $H_\infty$-calculus is an important tool in the 
modern analysis of nonlinear partial differential equations. The $H_\infty$-calculus 
was introduced by McIntosh in \cite{McIn}. We refer the reader to Denk, Hieber, Pr\"uss 
\cite{DHP} or Kunstmann, Weis \cite{KuWe} for recent and extensive surveys.  

Let $A:\calD(A)\subset Y\to Y$ be a closed and densely defined operator 
in a Banach space $Y$. In the sequel, $\Lambda=\Lambda(\theta)$ will denote a sector 
\begin{equation}\label{h2.B}
 \Lambda=\Lambda(\theta)=\{\lambda=re^{i\varphi}\st 
 r\ge 0,\;\theta\le\varphi\le 2\pi-\theta\},\qquad 0<\theta<\pi,
\end{equation}
in the complex plane. By $H=H(\theta)$ we denote the space of all holomorphic functions 
$\cz\setminus\Lambda\to\cz$ for which  
$|f(\lambda)|\le c(|\lambda|^\delta+|\lambda|^{-\delta})^{-1}$ for some $\delta>0$. 
If $\Lambda\setminus\{0\}$ is contained in the resolvent set of $A$ and   
$\|\lambda(\lambda-A)^{-1}\|_{\calL(Y)}$ is uniformly bounded on 
$\Lambda\setminus\{0\}$, then  
\begin{equation}\label{h2.A}
  f(A):=\frac{1}{2\pi i}\int_{\partial\Lambda} f(\lambda)(\lambda-A)^{-1}\,d\lambda,
  \qquad f\in H,
\end{equation}
converges absolutely and defines an element in $\calL(Y)$. By definition, the 
operator $A$ has a bounded 
$H_\infty$-calculus with respect to $\cz\setminus\Lambda$, if 
\begin{equation}\label{h2.C}
 \|f(A)\|_{\calL(Y)}\le M\,\|f\|_{\infty} \qquad\forall\;f\in H
\end{equation}
with a constant $M$ not depending on $f\in H$ (here, $\|f\|_\infty$ denotes the supremum 
norm of $f$). The name $H_\infty$-calculus originates from the fact that \eqref{h2.C} then 
allows the definition of $f(A)$ in $\calL(Y)$ for any bounded holomorphic function $f$ on 
$\cz\setminus\Lambda$. In particular, the choice of  $f(\lambda)=\lambda^{it}$, $t\in\rz$, 
implies the boundedness of purely imaginary powers, $A^{it}\in\calL(Y)$, and 
$\|A^{it}\|\le Me^{|t|\theta}$ for all $t\in\rz$.  
Establishing \eqref{h2.C} requires a thorough understanding of the resolvent of $A$ 
that goes well beyond proving the boundedness of $\|\lambda(\lambda-A)^{-1}\|_{\calL(Y)}$. 

In view of its importance for nonlinear parabolic equations, there is a vast literature 
concerned with the question of the existence of a bounded 
$H_\infty$-calculus for operators $A$ in quite different contexts. For example, 
Amann, Hieber, Simonett \cite{AHS} treat differential operators on $\rz^n$ and 
on compact manifolds with little regularity in the coefficients; 
Duong in \cite{Duon} considers boundary 
value problems on smooth manifolds, extending Seeley's work \cite{Seel2} on bounded 
imaginary powers; Denk, Dore, Hieber, Pr\"uss, and Venni \cite{DDHPV} then investigate 
boundary value problems of little regularity; in \cite{EsSe}, Escher and Seiler consider 
the Dirichlet-Neumann operator for domains of low regularity; \cite{CSS1} shows boundedness 
of the imaginary powers for differential operators on manifolds with conical singularities. 

In the present paper we study realizations (i.e.\ closed extensions) $\ulA_T$ of a 
$\mu$-th order cone differential operator $A$ subject to lower order differential 
boundary conditions $T$ on a manifold $\dz$ with conical singularities, where  
the boundary value problem $\spalte{A}{T}$ is assumed to be Shapiro-Lopatinskii 
elliptic. More precisely (for an explanation of the following notation see 
Section \ref{sectionh4}), we consider $A$ as an unbounded operator in a weighted 
$L_p$-space $\calH^{0,\gamma}_p(\dz,E)$ of sections of a vector bundle $E$, initially 
defined on a space of smooth sections that vanish under the boundary condition $T$. 
In general, $A$ will have a large number of 
closed extensions with domain contained in $\calH^{\mu,\gamma}_p(\dz,E)_T$, 
the space of sections of smoothness $\mu$ vanishing under $T$. 
We let $\ulA_T$ denote one of these extensions. 
All closed extensions can be described explicitly, cf.\ \cite{CSS2}. 

In Theorem \ref{hinfty} we show that $\ulA_T$ admits a bounded $H_\infty$-calculus 
provided the resolvent of $\ulA_T$ has the structure of an element of (a version of) the 
parameter-dependent cone calculus for pseudodifferential boundary value problems. 
In Section \ref{sectionh6} we next establish conditions on $\spalte{A}{T}$ which are 
more easily verified and guarantee that the resolvent of the minimal extension $A_{T,\min}$ 
has the desired structure. In essence, these conditions ask for the invertibility of 
all parameter-dependent principal symbols associated with $\spalte{\lambda-A}{T}$, for 
then the cone calculus developed by Schulze (see e.g.\ \cite{KaSc} for a presentation) 
allows the construction of a parametrix which yields the resolvent to $A_{T,\min}$. 
Based on the results on resolvents of closed extensions of 
cone differential operators in \cite{GKM2}, \cite{GKM1} by Gil, Krainer, and Mendoza, 
in \cite{Krai} by Krainer, 
and  in  \cite{ScSe2}, we expect to obtain 
such natural conditions also for extensions different from the minimal one. 
This analysis, however, is not the focus of the present paper. 

In Section \ref{sectionh7} we 
apply these results to the minimal extension of 
the Laplacian on a manifold with straight conical singularity, subject to Dirichlet or 
Neumann conditions. We find conditions on the dimension of $\dz$ and on the weight 
$\gamma$ such that $\Delta_{D,\min}$ and $\Delta_{N,\min}$ in
$\calH^{0,\gamma}_p(\dz)$ satisfy all the above assumptions. As a consequence, we derive 
maximal regularity for the associated initial boundary value problem. 

\forget{In their papers \cite{GKM2}, \cite{GKM1}, Gil, Krainer, and Mendoza study the resolvent 
of cone differential operators in the boundaryless case and, in particular, find 
conditions that imply the existence and decay property of the resolvent of an arbitrary 
closed extension in a sector $\Lambda$. Krainer in \cite{Krai} generalizes these results 
also to boundary value problems. }	

\section{Differential operators on smooth manifolds with boundary}\label{sectionh3}
In this section we want to recall classical results of \cite{Seel1}, \cite{Seel2} and 
\cite{Duon} on the resolvent of differential operators on manifolds with boundary and 
their $H_\infty$-calculus. 

In the following let $A$ be a differential operator of positive order $\mu\in\nz$ on 
a smooth manifold with boundary $X$, acting on sections into a vector bundle $E$,
 $$A:\;\ci(X,E)\longrightarrow \ci(X,E).$$
Moreover, let $T=(T_0,\ldots,T_{\mu-1})$ be a tuple of (normal) boundary conditions, 
i.e. 
 $$T_j=\gamma_0\circ B_j:\;\ci(X,E)\longrightarrow \ci(\partial X,F_j);$$
here, the $B_j$ are differential operators of order $\mu-j$ acting from sections 
into $E$ to sections into some hermitian vector bundle $F_j$ and $\gamma_0$ is 
the operator of restriction to the boundary of $X$. It is allowed that some of the 
$F_j$ are zero dimensional, i.e. the corresponding boundary condition $T_j$ is void. 

For a function space $\calF$ on $X$ let us set  
 $$\calF_T=\{u \in\calF\st Tu=0\},$$
provided application of $T$ to elements of $\calF$ makes sense. Let us now consider $A$ 
as the unbounded operator
\begin{equation}\label{h3.A} 
 A:\;\ci(X,E)_T\subset L_p(X,E)\longrightarrow L_p(X,E)
\end{equation}
with (a fixed) $1<p<\infty$. If the boundary value problem $\spalte{A}{T}$ 
is (Shapiro-Lopatinskii) elliptic, the closure $A_T$ of this 
operator can be shown to be defined by the action of $A$ on the domain 
\begin{equation}\label{h3.B} 
 \calD(A_T)=H^\mu_p(X,E)_T.
\end{equation}

\subsection{Elements of Boutet de Monvel's calculus for boundary value problems}
\label{sectionh3.1}

In \cite{Bout} Boutet de Monvel introduced a pseudodifferential calculus on manifolds 
with boundary which allowed the construction of parametrices to elliptic problems 
$\spalte{A}{T}$. This calculus has a corresponding parameter-dependent version, 
some of whose elements we describe now. For a short presentation see for example 
\cite{CSS2}. 

Let $\Sigma\subset\cz$ be a closed sector with vertex in the origin. We recall the 
definition of the space $B^{\nu,0}(X;\Sigma,E)$, $\nu\in\rz$, of operator-families 
on $X$, 
\begin{equation}\label{h3.C}
 A(\eta)=P_+(\eta)+G(\eta),
\end{equation}
depending on $\eta\in\Sigma$ as parameter. The first term is a pseudodifferential 
operator 
 $$P_+(\eta)=r_+\, P(\eta)\, e_+;$$
here, $e_+:\ci(X,E)\to\ci(2X,2E)$ denotes the operator of extending sections on $X$ 
by zero to sections on the double $2X$ (or any other smooth closed manifold containing 
$X$ as submanifold), while $r_+:\ci(2X,2E)\to\ci(X,E)$ is the operator of restriction. 
$P(\eta)$ is a classical, parameter-dependent pseudodifferential operator of order 
$\nu$ on $2X$. Furthermore it is required that $P(\eta)$ satisfies the 
{\em transmission condition} with respect to $\partial X$. Since differential operators 
as well as their parametrices always satisfy the transmission condition, we shall 
skip the details. 

The second term on the right-hand side of \eqref{h3.C} belongs to 
$B^{\nu,0}_G(X;\Sigma,E)$, the space of (singular) {\em Green operators}. 
Modulo {\em regularizing operators}, i.e.\ integral operators having a kernel in 
$\calS(\Sigma,\ci(X\times X))$, each Green operator is localized in a 
collar neighborhood of the boundary. By standard use of a partition of unity on $X$, this 
localized part is determined by operators on the half space 
$\overline{\rz}{}^n_+=\{(x^\prime,x_n)\st x^\prime\in\rz^{n-1},\,x_n\ge 0\}$
of the form 
 $$\op^\prime({d}\,)(\eta):\cicomp(\overline{\rz}{}^n_+,\cz^q)\longrightarrow
   \ci(\overline{\rz}{}^n_+,\cz^q),\qquad q=\text{dim}\,E,$$
defined by 
\begin{equation}\label{h3.D}
 \big(\op^\prime({d})(\eta)u\big)(x,t)=\int_{\rz^{n-1}}e^{ix\xi}\wt{d}(x,\xi,\eta;s,t)
 (\calF_{y\to\xi}u)(\xi,s)\,ds\dbar\xi
\end{equation}
with a {\em symbol kernel} 
\begin{equation}\label{h3.E}
 {d}(x,\xi,\eta;s,t)=\wt{d}(x,\xi,\eta;[\xi,\eta]s,[\xi,\eta]t),
\end{equation}
where $\wt{d}$ is a $(q\times q)$-matrix whose components are symbols 
\begin{equation}\label{h3.F}
 \wt{d}_{ij}(x,\xi,\eta;s,t)\in S^{\nu+1}_{\cl}(\rz^{n-1}\times\rz^{n-1}\times\Sigma,
 \calS(\rpbar\times\rpbar)).
\end{equation}
In \eqref{h3.E} and the sequel, $[\xi,\eta]$ denotes a smoothed norm function, i.e.\ 
a smooth positive function in $(\xi,\eta)$ that coincides with $|(\xi,\eta)|$ outside the 
unit ball. 
Roughly speaking, $\op^\prime({d})$ acts as an integral operator with smooth 
kernel in the direction normal to the boundary, while it is a pseudodifferential 
operator of order $\nu$ in the tangential direction. 

Both $B^{\nu,0}_G(X;\Sigma,E)$ and $B^{\nu,0}(X;\Sigma,E)$ are Fr\'{e}chet spaces in a 
natural way, since the spaces of local symbols and symbol kernels as well as the 
regularizing operators carry Fr\'{e}chet topologies.   

\subsection{Parameter-ellipticity and the resolvent of $A_T$}\label{sectionh3.2}

Let $\Lambda=\Lambda(\theta)$ be a sector in $\cz$ as introduced in \eqref{h2.B}. 
We shall call $A_T$, cf.\ \eqref{h3.A} and \eqref{h3.B}, 
(parameter-){\em elliptic with respect to} $\Lambda$ if two conditions are 
satisfied. First, the principal symbol of $A$ has no spectrum in $\Lambda$ for non-zero 
covariables, i.e. in local coordinates, 
\begin{equation}\label{h3.Z}
 \mbox{\rm det}\,(\lambda-\sigma^\mu_\psi(A)(x,\xi))\not=0\qquad\forall\;
 0\not=(\xi,\lambda)\in\rz^n_{\xi}\times\Lambda.
\end{equation}
The second condition is the invertibility of the boundary symbol of 
$\spalte{\lambda-A}{T}$; in local coordinates $(x^\prime,x_n)$ near the boundary of 
$X$ this means that, for fixed $(x^\prime,\xi^\prime,\lambda)$, 
 $$\begin{pmatrix}
    \lambda-\sigma^\mu_\psi(A)(x^\prime,0,\xi^\prime,D_{x_n})\\ 
    \left(\gamma_0\circ\sigma^{j}_\psi(B_j)
    (x^\prime,0,\xi^\prime,D_{x_n})\right)_{j=0,\ldots,\mu-1}
   \end{pmatrix}:\;
   \calS(\rpbar,\cz^q)\longrightarrow
   \begin{matrix}
    \calS(\rpbar,\cz^q)\\
    \oplus\\
    \cz^{m_0}\oplus\ldots\oplus\cz^{m_{\mu-1}}
   \end{matrix}$$
is an isomorphism whenever $(\xi^\prime,\lambda)\not=0$. Here, 
$m_j=\mbox{\rm dim}\,F_j$ and $\gamma_0$ denotes evaluation at 0 of functions on 
$\rpbar$. Now let 
\begin{equation}\label{h3.G} 
 \Sigma=\Sigma(\theta,\mu)=\left\{\eta=se^{i\alpha}\st s\ge 0,\;
 \mbox{$\frac{\theta}{\mu}\le\alpha\le\frac{2\pi-\theta}{\mu}$}\right\}. 
\end{equation}
Then $\Sigma$ is a sector and $\Lambda=\{\eta^\mu\st\eta\in\Sigma\}$ for $\Lambda$ 
as in \eqref{h2.B}. 

\begin{theorem}\label{h3.1}
Let $A_T$ be elliptic with respect to $\Lambda$. Then 
 $$\lambda-A_T:\,H^\mu_p(X,E)_T\subset L_p(X,E)\longrightarrow L_p(X,E)$$
is invertible for large $\lambda\in\Lambda$ and there exists an operator-family 
$P_+(\eta)+G(\eta)\in B^{-\mu,0}(X;\Sigma,E)$ such that, for large $\eta\in\Sigma$,
 $$(\eta^\mu-A_T)^{-1}=P_+(\eta)+G(\eta).$$
\end{theorem}

This theorem was essentially proven by Seeley \cite{Seel1} but without 
using the `language' of Boutet de Monvel's calculus. For other proofs we 
refer to \cite{Grub} and \cite{Grub2}. 

\begin{corollary}\label{h3.2}
If $A_T$ is elliptic with respect to $\Lambda$ then there exists a constant $c_p\ge0$ 
such that for all large $\lambda\in\Lambda$ 
 $$\|(\lambda-A_T)^{-1}\|_{L_p(X,E)}\le c_p\,\frac{1}{|\lambda|}.$$
\end{corollary}

Using the above resolvent structure, \cite{Seel2} shows the existence of bounded 
imaginary powers, while \cite{Duon} proves a bounded $H_\infty$-calculus:  

\begin{theorem}\label{h3.3}
If $A_T$ is elliptic with respect to $\Lambda$ then there exists a constant $c\ge0$ 
such that $A_T+c$ has a bounded $H_\infty$-calculus with respect to $\cz\setminus\Lambda$ 
$($simultaneously for all $1<p<\infty)$. 
\end{theorem}

\section{Boundary value problems on conic manifolds}\label{sectionh4}
The main objective of this section is the description of boundary value problems on a 
manifold with conical singularities and of a parameter-dependent Boutet de Monvel 
calculus adapted to this situation. For simplicity we shall confine the description 
to the case of one conical singularity. 

\subsection{Cone differential operators}\label{sectionh4.0}

Let $\intd$ be an $(n+1)$-dimensional riemannian manifold with boundary having a
cylindrical end, i.e.\ there exists a compact subset $C$ such that $\intd\setminus C$ is 
isometric to the product ${]0,1[}\times X$ for a smooth compact (not necessarily 
connected) riemannian manifold $X$ with boundary. 
We fix the coordinate in ${]0,1[}$ in such a way that every neighborhood 
of $C$ in $\intd$ has nonempty intersection with ${]\frac{1}{2},1[}\times X$. 
We complete $\intd$ with the help of the riemannian metric. If $\dz$ is the resulting space, 
then $\dz\setminus C$ can be identified with ${[0,1[}\times X$. We next denote by 
$\intb$ the boundary of $\intd$ and by $\bz$ its completion in the metric 
inherited from $\intd$. Then $\bz\setminus C$ can be identified with 
${[0,1[}\times\partial X$ and $\bz$ itself is a smooth manifold with boundary.  

\forget{
To make the exposition simpler we shall assume in the following that $\dz$ has only one 
cylindrical end, denoted by ${[0,1[}\times X$. Here we use the canonical coordinates $(t,x)$, 
$0\le t<1$, $x\in X$. The subset $\{0\}\times X$ of $\dz$, where $t=0$, is  occasionally 
called the {\em singularity} of $\dz$. 
}

In the sequel, a vector bundle over $\dz$ will be a smooth hermitian vector bundle over 
$\intd$ such that $E|_{{]0,1[}\times X}$ is isometric to the pull-back under the 
canonical projection ${]0,1[}\times X\to X$ of a hermitian bundle $E_0$ over $X$. 

That we call $\dz$ a {\em manifold with conical singularity} is due to the class of 
differential operators we consider on it. A $\mu$-th order differential operator $A$ on 
$\intd$ with smooth coefficients, acting between sections of $E$, is called a 
{\em cone differential operator} if   
\begin{equation}\label{fuchs}
 A=t^{-\mu}\smsum_{j=0}^\mu a_j(t)(-t\partial_t)^j,\qquad
 a_j\in\ci([0,1[,\text{\rm Diff}^{\,\mu-j}(X;E_0)),
\end{equation}
on ${]0,1[}\times X$; here we use the canonical coordinates $(t,x)$ with 
$0\le t<1$, $x\in X$, and $\text{\rm Diff}^{k}(X;E_0)$ denotes the space of $k$-th order  
differential operators on $X$ with smooth coefficients, acting between sections of $E_0$.  

\subsection{Function spaces}\label{sectionh4.1}

We first introduce some function spaces on $\dz$ (the corresponding spaces on $\bz$ 
are obtained analogously). Patching together two copies of 
$\dz$, we obtain a smooth manifold with boundary $2\dz$ containing $\dz$. 
Let $\omega$ be an arbitrary fixed cut-off function on $\dz$ (extended by 1 to $2\dz$). 

\begin{definition}\label{h4.1}
Let $\gamma\in\rz$. Then $\cig(\dz)$ denotes the space of all smooth functions $u$ 
for which $(1-\omega)u\in\ci(2\dz)$ and 
\begin{equation}\label{h4.X}
 t\mapsto t^{\frac{n+1}{2}-\gamma}(\log t)^k (t\partial_t)^l (\omega u)(t,\cdot)
\end{equation}
is a bounded function on ${]0,1[}$ with values in $\ci(X)$ for all $k,l\in\nz_0$. 
Moreover, we set 
 $$\cii(\dz)=\schnitt_{\gamma\in\rz}\cig(\dz).$$
\end{definition}

\begin{definition}\label{h4.2}
For $s\in\nz_0$, $\gamma\in\rz$, and $1<p<\infty$, we let $\calH^{s,\gamma}_p(\dz)$ 
be the space of all functions $u$ such that $(1-\omega)u\in H^s_p(2\dz)$ and 
 $$t^{\frac{n+1}{2}-\gamma}(t\partial_t)^l\partial^\alpha_x (\omega u)(t,x)\;\in\;
   L_p({]0,1[}\times X,\mbox{$\frac{dt}{t}dx$})\qquad\forall\;l+|\alpha|\le s,$$
where $dx$ refers to some metric on the manifold with boundary $X$.  
\end{definition}

These definitions naturally extend to sections into $E$ and to real numbers $s$ 
(however, we shall only need the case of nonnegative integers). Then we write 
$\cig(\dz,E)$, $\cii(\dz,E)$, and $\calH^{s,\gamma}_p(\dz,E)$. 

\subsection{Elliptic boundary value problems}\label{sectionh4.2}

Let $A$ be a $\mu$-th order cone differential operator on $\dz$ as described in the 
introduction, 
 $$A:\;\cii(\dz,E)\longrightarrow\cii(\dz,E),$$
and $T=(T_0,\ldots,T_{\mu-1})$ be a tuple of (normal) boundary conditions, i.e. 
\begin{equation}\label{h4.A.6}
 T_j=\gamma_0\circ B_j:\;\cii(\dz,E)\longrightarrow \cii(\bz,F_j)
\end{equation}
with cone differential operators $B_j$ of order $j$, and $\gamma_0$ being the operator of 
restriction to the boundary $\bz$ of $\dz$. Setting $F=F_0\oplus\ldots\oplus F_{\mu-1}$ we 
thus have 
\begin{equation}\label{h4.A.2}
 \begin{pmatrix}A\\T\end{pmatrix}:\; 
 \cii(\dz,E)\longrightarrow
 \begin{matrix}\cii(\dz,E)\\\oplus\\ \cii(\bz,F)\end{matrix}.
\end{equation}
Let us shortly sketch what $\dz$-ellipticity of $\spalte{A}{T}$ means; 
for details we refer to 
Section 3.2 of \cite{CSS2}. First, $A$ is elliptic on $\intd$ in the standard sense, 
i.e.\ the  principal symbol $\sigma^\mu_\psi(A)$ is everywhere 
(i.e.\ for non-zero covariables) invertible. 
In the splitting of coordinates $(t,x)\in{]0,1]}\times X$ near $t=0$ we have 
 $$\sigma^\mu_\psi(A)(t,x,\tau,\xi)=t^{-\mu}\smsum_{j=0}^\mu
   \sigma^{\mu-j}_\psi(a_j)(t,x,\xi)(-it\tau)^j.$$
Observe a characteristic ``degenerate'' structure: There is the singular factor $t^{-\mu}$ 
and the covariable $\tau$ appears only in the form $t\tau$. Removing this degeneracy and 
freezing the coefficients in $t=0$ we obtain the so-called {\em rescaled symbol} 
 $$\wt\sigma^\mu_\psi(A)(x,\tau,\xi)=\smsum_{j=0}^\mu
   \sigma^{\mu-j}_\psi(a_j)(0,x,\xi)(-i\tau)^j.$$
We require this symbol to be invertible for all $x$ and $(\tau,\xi)\not=0$. 
Since $\spalte{A}{T}$ is a usual boundary value problem on $\intd$, we may associate 
with it the standard principal boundary symbol. Ellipticity asks the invertibility of 
this symbol. Since the boundary conditions are also given by cone differential operators, 
the boundary symbol again has a degenerate structure. Removing the degeneracy and 
freezing coefficients in $t=0$ yields the {\em rescaled boundary symbol} which is also 
required to be invertible. 

\subsection{Realizations of $A$ with respect to $T$}\label{sectionh4.3}

For a function space $\calF$ on $\dz$ we use again the notation 
$\calF_T=\{u\in\calF\st Tu=0\}$. Let us now consider the unbounded operator 
\begin{equation}\label{h4.A}
 A:\;\cii(\dz,E)_T\subset\calH^{0,\gamma}_p(\dz,E)\longrightarrow
 \calH^{0,\gamma}_p(\dz,E)
\end{equation}
with $\gamma\in\rz$ and $1<p<\infty$ (it would be more precise to write 
$A_{\gamma,p}$ but for convenience we exclude $\gamma$ and $p$ from the notation). 
We shall assume that $\spalte{A}{T}$ is $\dz$-elliptic in the sense of 
Section \ref{sectionh4.2}. 

We write $A_{T,\min}$ for the closure of $A$ and define $A_{T,\max}$ by the action of 
$A$ on the space 
 $$\calD(A_{T,\max})=\{u\in\calH^{\mu,\gamma}_p(\dz,E)_T\st 
   Au\in \calH^{0,\gamma}_p(\dz,E)\}.$$
By abuse of notation, $A_{T,\max}$ is not the true maximal extension of $A$, but it is 
the largest one that still `feels' the boundary condition $T$. The closed extensions 
between the minimal and the maximal are usually called the {\em realizations} of $A$ 
with respect to $T$. In \cite{CSS2}, Theorem 5.12, we have shown: 

\begin{theorem}\label{h4.3.5} 
There exists a finite dimensional space
$\calE=\calE_{A,T}\subset\calC^{\infty,\gamma}(\dz,E)$ 
of smooth functions on $\intd$ such that  
 $$\calD(A_{T,\max})=\calD(A_{T,\min})\oplus\calE.$$
\end{theorem}

As a consequence, any realization $\ulA_T$ is determined by a subspace 
$\underline{\calE}$ of $\calE$, i.e. 
 $$\calD(\ulA_T)=\calD(A_{T,\min})\oplus\underline{\calE}.$$
In \cite{CSS2} we also have characterized the domain of the closure, namely 
 $$\calD(A_{T,\min})=\Big\{u\in\schnitt_{\eps>0}\calH^{\mu,\gamma+\mu-\eps}_p(\dz,E)_T\st 
   Au\in\calH^{0,\gamma}_p(\dz,E)\Big\}.$$
In particular, $\calH^{\mu,\gamma+\mu}_p(\dz,E)_T\subset\calD(A_{T,\min})
\subset\calH^{\mu,\gamma+\mu-\eps}_p(\dz,E)_T$ for any $\eps>0$. If in addition $A_T$ is 
{\em conormally elliptic} with respect to the weight $\gamma+\mu$, see (E3) in 
Section \ref{sectionh6} for an explanation, we even have 
\begin{equation}\label{h4.A.5}
 \calD(A_{T,\min})=\calH^{\mu,\gamma+\mu}_p(\dz,E)_T.
\end{equation}

\subsection{Parameter-dependent operators}\label{sectionh4.4}

A generalization of Boutet de Monvel's calculus to boundary value problems on 
manifolds with conical singularities was introduced in \cite{ScSc1} and \cite{ScSc2}. 
A corresponding parameter-dependent version can be found e.g.\ in \cite{KaSc}. 
We give here a somewhat simplified presentation of some of the elements of this 
calculus, focusing on the structures that are necessary for the description 
of resolvents. For convenience, we shall assume $E=\cz$, since the general case is 
a straightforward extension of this situation.  

\subsubsection{Green operators}\label{sectionh4.4.1}

With $X^\wedge:=\rz_+\times X$ and for $\gamma\in\rz$ let us set
\begin{align} 
 \calH^{0,\gamma}_p(X^\wedge)&=
 L_p(X^\wedge,t^{(\frac{n+1}{2}-\gamma)p}\mbox{$\frac{dt}{t}$}dx)=
   t^{\gamma-\frac{n+1}{2}} L_p(X^\wedge,\mbox{$\frac{dt}{t}$}dx),\label{h4.B.2}\\
 \calS^\gamma_0(X^\wedge)&=\{u\in\ci(X^\wedge)\st\omega u\in\calC^{\infty,\gamma}(\dz)
   \text{ and }(1-\omega)u\in\calS(\rz,\ci(X))\}.\label{h4.B.4}
\end{align}
The latter is a Fr\'echet space in a natural way. 

Let $\Sigma$ be a closed subsector of $\cz$ with vertex in 0 and $\nu,\gamma\in\rz$. 
Then we denote by $R^{\nu,0}_G(X^\wedge;\Sigma,\gamma)$ the space of all 
operator families $G(\eta)$, $\eta\in\Sigma$, such that 
 $$\big(G(\eta)u\big)(t,x)=
   \big\langle k(\eta;t,x),\overline{u}\big\rangle_{\calH^{0,0}_2(X^\wedge)}=
   \int_{X^\wedge}k(\eta;t,x,t^\prime,x^\prime)u(t^\prime,x^\prime)\,
   t^{\prime n}dt^\prime dx^\prime,$$
where the integral kernel $k$ has the form 
\begin{equation}\label{h4.B}
 k(\eta;t,x,t^\prime,x^\prime)=\wt{k}(\eta;[\eta]t,x,[\eta]t^\prime,x^\prime)
\end{equation}
with 
\begin{equation}\label{h4.C}
 \wt{k}(\eta;t,x,t^\prime,x^\prime)\in S^{\nu+n+1}\big(\Sigma,
 \calS^{\gamma+\eps}_0(X^\wedge_{(t,x)})\pit
 \calS^{-\gamma+\eps}_0(X^\wedge_{(t^\prime,x^\prime)})\big)
\end{equation}
for some $\eps>0$ (depending on $G$). 

\subsubsection{Mellin pseudodifferential operators}\label{sectionh4.4.2}

For $\beta\in\rz$ let $\Gamma_\beta=\{z\in\cz\st\re z=\beta\}$. We denote by 
$M^{\nu,0}_\beta(X;\Sigma)$, $\nu\in\rz$, the space of all functions 
$h:\Gamma_{\beta}\times\Sigma\to B^{\nu,0}(X)$ such that 
 $$h(\mbox{$\beta+i\tau$},\eta)\in B^{\nu,0}(X;\rz_\tau\times\Sigma_\eta)$$
(cf.\ Section \ref{sectionh3.1} with parameter space $\rz\times\Sigma$ instead of 
$\Sigma$). A {\em Mellin symbol} is a function 
$h=h(t,\beta+i\tau,\eta)\in\ci(\rz_+,M^{\nu,0}_\beta(X;\Sigma))$; it induces the 
Mellin pseudodifferential operator 
$\op_M^{\frac{1}{2}-\beta}(h)(\eta):\;\calC^\infty_{\text{\rm comp}}(X^\wedge)
\to\calC^\infty(X^\wedge)$
by 
 $$\big(\op_M^{\frac{1}{2}-\beta}(h)(\eta)u\big)(t)=
   \int_{-\infty}^\infty t^{-\beta-i\tau}h(t,\beta+i\tau,\eta)
   (\calM u)(\beta+i\tau)\,\dbar\tau;$$
here we have identified $\calC^\infty(X^\wedge)$ with $\calC^\infty(\rz_+,\ci(X))$, 
and $\calM$ denotes the Mellin transform, i.e. 
$(\calM u)(z)={\displaystyle\int_0^\infty} t^z u(t)\,\frac{dt}{t}$ for $z\in\cz$. 

\subsubsection{The full class and a norm estimate}\label{sectionh4.4.3}

We define ${C}^{\nu,0}(\dz;\Sigma,\gamma)$ with $\nu\le 0$ and $\gamma\in\rz$ as the set 
of all operator-families $A(\eta)$ of the form 
\begin{equation}\label{h4.D}
 A(\eta)=\sigma_0\left\{
           t^{-\nu}\op_M^{\gamma-\frac{n}{2}}(h)(\eta)+G(\eta)
           \right\}\sigma_1+
           (1-\sigma_3)P(\eta)(1-\sigma_4)+G_\infty(\eta),
\end{equation}
where $\sigma_j\in\cicomp([0,1[)$ are cut-off functions and 
\begin{itemize}
 \item[a)] $h(t,z,\eta)=\wt{h}(t,z,t\eta)$ with 
  $\wt{h}\in \ci(\rpbar,M^{\nu,0}_{\frac{n+1}{2}-\gamma}(X;\Sigma))$, 
  cf.\ Section \ref{sectionh4.4.2},
 \item[b)] $P(\eta)\in B^{\nu,0}(2\dz;\Sigma)$, cf.\ Section \ref{sectionh3.1}, 
 \item[c)] $G(\eta)\in R^{\nu,0}_G(X^\wedge;\Sigma,\gamma)$, 
  cf.\ Section \ref{sectionh4.4.1},
 \item[d)] $G_\infty\in\calS(\Sigma,\calL(\calH^{s,\gamma}_p(\dz),
  \calH^{s^\prime,\gamma}_p(\dz)))$ for all $1<p<\infty$ and $s,s^\prime\in\rz$ 
  with $s>-1+\frac{1}{p}$. 
\end{itemize}
Initially, each such operator is defined on smooth functions, i.e.\ 
$A(\eta):\cicomp(\intd)\to\ci(\intd)$, where multiplications by $\sigma_1$ and 
$\sigma_0$ are viewed as operators $\cicomp(\intd)\to\cicomp(X^\wedge)$ and 
$\cicomp(X^\wedge)\to\cicomp(\intd)$, respectively. Similarly, $1-\sigma_3$ and 
$1-\sigma_4$ act between smooth functions on $\dz$ and $2\dz$. 

Let $A(\eta)\in{C}^{\nu,0}(\dz;\Sigma,\gamma)$ as in \eqref{h4.D}. 
We deduce from \eqref{h3.C} that there is a decomposition of $\wt{h}$ as     
\begin{equation}\label{h4.D.5}
 \wt{h}(t,\mbox{$\frac{n+1}{2}$}-\gamma+i\tau,\eta)=
 \wt{p}_+(t,\tau,\eta)+\wt{g}(t,\tau,\eta)
\end{equation}
where $t\mapsto\wt{p}(t,\tau,\eta)$ is a smooth (up to $t=0$) family of 
parameter-dependent pseudodifferential operators of order $\nu$ on $2X$ with parameter 
space $\rz\times\Sigma$, while $\wt{g}(t,\tau,\eta)$ belongs to 
$\ci(\rpbar,B^{\nu,0}_G(X;\rz\times\Sigma))$. In local coordinates, 
\begin{equation}\label{h4.E}
 \wt{p}_+(t,\tau,\eta)=\op_+(\wt{a})(t,\tau,\eta)\text{ with }
 \wt{a}(t,x,\tau,\xi,\eta)\in 
 S^{\nu}_\cl(\rpbar\times\rz^n_x\times\rz^{1+n}_{(\tau,\xi)}\times\Sigma)
\end{equation}
and\footnote{By definition, 
$S^{\nu}_\cl(\rpbar\times\rz^n_x\times\rz^{1+n}_{(\tau,\xi)}\times\Sigma)$ denotes the 
space of all functions $b$ that are smooth on 
$\rpbar\times\rz^n_x\times\rz^{1+n}_{(\tau,\xi)}\times\Sigma$ satisfy there uniform estimates 
$|\partial_t^l\partial_x^\beta\partial^\alpha_{(\tau,\xi,\eta)}b(t,x,\tau,\xi,\eta)|\le
c\spk{\tau,\xi,\eta}^{\nu-|\alpha|}$ for any order of derivatives and that have 
asymptotic expansions into components that are positively homogeneous in 
$(\tau,\xi,\eta)$.}  
$\wt{g}(t,\tau,\eta)$ is defined via a symbol kernel 
\begin{equation}\label{h4.F}
 \wt{d}(t,x^\prime,\tau,\xi^\prime,\eta;u,v)=
 d(t,x^\prime,\tau,\xi^\prime,\eta;[\tau,\xi^\prime,\eta]u,[\tau,\xi^\prime,\eta]v)
\end{equation}
with 
\begin{equation}\label{h4.G}
 d(t,x^\prime,\tau,\xi^\prime,\eta;u,v)\in
 S^{\nu+1}(\rpbar\times\rz^{n-1}_{x^\prime}\times\rz^{n}_{(\tau,\xi^\prime)}\times\Sigma;
 \calS(\rpbar\times\rpbar)).
\end{equation}
The non-local, smoothing terms arising are as above in d). 
 
\begin{theorem}\label{h4.2.4}
Let $A(\eta)\in{C}^{\nu,0}(\dz;\Sigma,\gamma)$, $\nu\le0$. Then $A(\eta)$ extends 
for each $\eta$ to a bounded operator in $\calH^{0,\gamma}_p(\dz)$ and there exists 
a constant $c_p\ge0$ such that 
 $$\|A(\eta)\|_{\calL(\calH^{0,\gamma}_p(\dz))}\le c_p\,\spk{\eta}^\nu\qquad
   \forall\;\eta\in\Sigma.$$
\end{theorem}
\begin{proof}
The proof is a combination of that for the boundaryless case (cf.\ Proposition 1 in 
\cite{CSS1}) and the results of \cite{GrKo} on the $L_p$-continuity of operators from 
Boutet's calculus. To give an idea, we shall derive the desired norm estimate for the 
operator $\sigma_0\,t^{-\nu}\,\op_M^{\gamma-\frac{n}{2}}(g)\,\sigma_1$ with 
$g(t,\tau,\eta)=\wt{g}(\tau,t\eta)$ as in \eqref{h4.D.5} (the $t$-independence of 
$\wt{g}$ can always be reached by a tensor product argument). By conjugation with 
$t$-powers, it is no restriction to assume that $\gamma=(n+1)(\frac{1}{2}-\frac{1}{p})$. 
With this choice of $\gamma$, we have 
$\calH^{0,\gamma}_p(X^\wedge)=L_p(X^\wedge,t^ndtdx)$ and the operators 
 $$\kappa_{\varrho}\in\calL(\calH^{0,\gamma}_p(X^\wedge)),\text{ defined by } 
   (\kappa_{\varrho}u)(t,x)=\varrho^{\frac{n+1}{p}}u(\varrho t,x),$$
provide bijective isometries for each $\varrho>0$. Moreover, defining 
$(Su)(r,x)=e^{(\frac{n}{2}-\gamma)r}u(e^{-r},x)$ for functions $u$ on 
$X^\wedge$, we obtain an isometric isomorphism 
 $$S:\calH^{0,\gamma}_p(X^\wedge)\longrightarrow L_p(\rz\times X,drdx).$$
Thus we are done if we can prove that 
 $$\spk{\eta}^{-\nu}\|S\,\kappa_{\spk{\eta}}^{-1}\,t^{-\nu}\,
   \op_M^{\gamma-\frac{n}{2}}(g)(\eta)\,\kappa_{\spk{\eta}}\,
   S^{-1}\|_{\calL(L_p(\rz\times X,drdx))}\le c_p, $$
uniformly in $\eta$. Now note that conjugation of $\op_M^{\gamma-\frac{n}{2}}(g)(\eta)$ 
by $\kappa_{\spk{\eta}}$ amounts to replacing ${g}(t,z,\eta)$ by 
${g}(\frac{t}{\spk{\eta}},z,\eta)$, and conjugation with $S$ transforms a Mellin 
operator on $X^\wedge$ with symbol $h(t,\frac{n+1}{2}-\gamma+i\tau)$ into a  
pseudodifferential operator on $\rz\times X$ with symbol 
$a(r,\varrho)=h(e^{-r},\frac{n+1}{2}-\gamma+i\varrho)$. 
Passing to the local situation, cf.\ \eqref{h4.F} and \eqref{h4.G}, we have to show that 
$\|\op^\prime(d_0)(\eta)\|_{\calL(L_p(\overline{\rz}^{n+1}_+))}\le c_p$ 
for the symbol kernel 
 $$d_0(r,x^\prime,\tau,\xi^\prime,\eta;u,v)=e^{r\nu}
   \mbox{$\wt{d}(x^\prime,\tau,\xi^\prime,e^{-r}\frac{\eta}{\spk{\eta}};
   [\tau,\xi^\prime,e^{-r}\frac{\eta}{\spk{\eta}}]u,
   [\tau,\xi^\prime,e^{-r}\frac{\eta}{\spk{\eta}}]v)$},$$
see \eqref{h3.D} for the definition of $\op^\prime$. By a straightforward calculation,  
 $$\|u^{k^\prime}D^k_u v^{l^\prime}D^l_v D^{m^\prime}_r D^m_\tau 
   D^{\alpha^\prime}_{x^\prime} D^\alpha_{\xi^\prime}
   d_0(r,x^\prime,\tau,\xi^\prime,\eta;u,v)\|_{L^2(\rz_{+,u}\times\rz_{+,v})}\le 
   C\,\spk{\tau,\xi^\prime}^{-|\alpha|-m+k-k^\prime+l-l^\prime},$$
uniformly in $(r,x^\prime,\tau,\xi^\prime,\eta)$. Now Theorem 4.1.(5) in \cite{GrKo} 
gives the continuity of $\op^\prime(d_0)(\eta)$ in $L_p(\overline{\rz}^{n+1}_+)$, with 
operator norm uniformly bounded in $\eta$. 
\end{proof}

The local symbols $\wt{a}$ in \eqref{h4.E} are defined for 
$(\tau,\xi,\eta)\in\rz^{1+n}\times\Sigma$. The symbols that arise in the analysis 
of resolvents extend -- holomorphically in $\eta$ -- to larger subsets of $\cz$. 
We shall need this property in particular for the principal part and define a 
corresponding class: 

\begin{definition}\label{h4.3}
For $c_0>0$ let $\Omega_{(\tau,\xi)}=\{\eta\in\cz\st|\eta|\ge c_0\spk{\tau,\xi}\}$. 
The set $C_{\rm hol}^{\nu,0}(\dz;\Sigma,\gamma)$ consists of all operator-families 
$A(\eta)\in{C}^{\nu,0}(\dz;\Sigma,\gamma)$ such that the local symbols $\wt{a}$ from 
\eqref{h4.E} admit a decomposition $\wt{a}=\wt{a}_0+\wt{a}_1$ with the following 
properties: 
\begin{itemize}
 \item[a)] $\wt{a}_0$ extends holomorphically to $\eta\in\Omega_{(\tau,\xi)}$ 
  $($for some $c_0)$ and satisfies the estimates 
  \begin{equation}\label{h4.H}
   |\partial_t^l\partial_x^\beta\partial_\tau^k\partial_\xi^\alpha\partial_\eta^\gamma\,
   \wt{a}_0(t,x,\tau,\xi,\eta)|\le 
   c\spk{\tau,\xi,\eta}^{\nu-k-|\alpha|-|\gamma|}
  \end{equation}
  uniformly in $(\tau,\xi)\in\rz^{1+n}$ and $\eta\in\Omega_{(\tau,\xi)}\cup\Sigma$ 
  for any order of derivatives, 
 \item[b)] $\wt{a}_1(t,x,\tau,\xi,\eta)\in 
   S^{\nu-1}_\cl(\rpbar\times\rz^n_x\times\rz^{1+n}_{(\tau,\xi)}\times\Sigma_\eta)$. 
\end{itemize}
For the pseudodifferential symbols of $P(\eta)$ in \eqref{h4.D} we require the 
corresponding structure.  
\end{definition}

\section{$H_\infty$-calculus of cone differential operators}\label{sectionh5}
The aim of this section is to prove the following theorem:  

\begin{theorem}\label{hinfty}
If the resolvent set of the realization $\ulA_T$ contains 
$\{\lambda\in\Lambda\st |\lambda|\ge R\}$ for some $R\ge0$ and there exists an 
$A(\eta)\in C^{-\mu,0}(\dz;\Sigma,\gamma)$ such that 
 $$(\eta^\mu-\ulA_T)^{-1}=A(\eta)\qquad\forall\;|\eta|\ge R^{1/\mu},$$
then there exists a constant $c\ge0$ such 
that $c+\ulA_T$ admits a bounded $H_\infty$-calculus with respect to $\cz\setminus\Lambda$ 
$($simultaneously for all $1<p<\infty)$. 
\end{theorem}

Before going into details, let us give an outline of the proof. In large parts it 
follows the proof of Theorem 5.1 in \cite{CSS1}, where we showed the existence of 
bounded imaginary powers for operators on conic manifolds without boundary. 

By replacing from the very beginning the differential 
operator $A$ by $c+A$ we can assume that the resolvent set contains all of $\Lambda$ and 
 $$A(\eta):=(\eta^\mu-\ulA_T)^{-1}\in C^{-\mu,0}(\dz;\Sigma,\gamma).$$ 
We then show that $\ulA_T$ itself has a bounded $H_\infty$-calculus. 
Let us remark that, for $f\in H$, 
 $$f(\ulA_T)=\int_{\partial\Lambda}f(\lambda)(\lambda-\ulA_T)^{-1}\,d\lambda=
   \int_{\partial\Sigma}f(\eta^\mu)A(\eta)\eta^{\mu-1}\,d\eta.$$
Inserting the explicit formula for $A(\eta)$, cf.\ \eqref{h4.D}, we obtain four 
integrals over the boundary of $\Sigma$ which have to be estimated from above by 
$M\|f\|_\infty$ with a constant $M$ independent of $f\in H$, cf.\eqref{h2.C}. 

Obviously, the integral associated with $G_\infty(\eta)$ can be estimated as desired, 
since $\eta^{\mu-1}G_\infty(\eta)$ is reapidly decreasing in $\eta$, hence integrable. 

For the integral associated with $(1-\sigma)P(\eta)(1-\sigma_1)$ we can apply 
Theorem \ref{h3.3}, since away from the singularity the Sobolev spaces and 
operator classes coincide with the usual ones on the smooth manifold with 
boundary $2\dz$. 

In Propositions \ref{lower}, \ref{principal}, and \ref{singg} we shall treat the term   
\begin{equation}\label{mellinterm}
 \sigma_0\int_{\partial\Sigma} f(\eta^\mu) \, t^\mu\opm{\gamma-\frac{n}{2}}(h)(\eta)
 \,\eta^{\mu-1}\,d\eta\,\sigma_1=
 \opm{\gamma-\frac{n}{2}}\Big(\sigma_0\,t^\mu\int_{\partial\Sigma} f(\eta^\mu)\,h(\eta)
 \,\eta^{\mu-1}\,d\eta\Big)\sigma_1.
\end{equation} 
First, however, we shall study the term induced by $G(\eta)$. 
To this end note that  multiplication by a cut-off function 
$\sigma\in\cicomp([0,1[)$ induces continuous 
operators $\calH^{0,\gamma}_p(\dz)\to\calH^{0,\gamma}_p(X^\wedge)$ and 
$\calH^{0,\gamma}_p(X^\wedge)\to\calH^{0,\gamma}_p(\dz)$, cf.\ \eqref{h4.B.2}. 
Thus it suffices to show the following:  


\begin{proposition}\label{greenpart}
Let $G(\eta)\in R^{-\mu,0}_G(X^\wedge;\Sigma,\gamma)$ and 
$G_f=\sigma_0 \displaystyle\int_{\partial\Sigma}
f(\eta^\mu) G(\eta)\eta^{\mu-1}\,d\eta \, \sigma_1$ for 
$f\in H$. Then $G_f\in\calL(\calH^{0,\gamma}_p(X^\wedge))$, and there 
exists a constant $M_p\ge0$ such that  
\begin{equation}\label{normbound}
 \|G_f\|_{\calL(\calH^{0,\gamma}_p(X^\wedge))}\le 
 M_p\| f \|_{\infty}\qquad\forall\; f\in H. 
\end{equation}
\end{proposition}
\begin{proof}
By conjugation with $t^\gamma$ we can assume that $\gamma=0$.  
Due to the symmetry of $\partial\Sigma$ and the fact that $\eta^{\mu-1}G(\eta)$ is 
integrable on compact parts of $\partial\Sigma$, we may confine ourselves to the 
integration over the set  
  $$\calC=\{\varrho e^{i\alpha}\st 1\le\varrho<\infty\}.$$
According to Section \ref{sectionh4.4.1}, $G(\eta)$ is 
an integral operator (with respect to the scalar product in 
$\calH^{0,0}_2(X^\wedge)$). Suppressing, for notational simplicity, 
the dependence on the variables $x$, $x^\prime$, the kernel of 
$G(\eta)$, $|\eta|\ge1$ is $k(\eta,t,s)=\wt{k}(\eta,|\eta|\,t,|\eta|\,s)$,  
where, for some $\varepsilon>0$, 
  $$\wt{k}(\eta,t,s)\in S^{-\mu+n+1}(\Sigma) \pit 
    \calS^\eps_0(X^\wedge)\pit \calS^\eps_0(X^\wedge).$$
Then $G_f$ the an integral operator with kernel 
\begin{equation}\label{kernelgz}
   k_f(t,s)=\sigma_0(t)\,\sigma_1(s)\int_\calC 
   f(\eta^\mu)k(\eta,t,s)\,\eta^{\mu-1} \, d\eta.
\end{equation} 
Writing 
$\wt{k}(\eta,t,s)=({\chi}(t)+(1-{\chi})(t))\wt{k}(\eta,t,s)({\chi}(s)+(1-{\chi})(s))$ 
with the characteristic function ${\chi}$ of $[0,1]$, the proposition will be 
true, if we can show that in each one of the four cases 
 \begin{align}
  k_f(t,s)&=\sigma_0(t)\,\sigma_1(s)\int_{\calC} f(\eta^\mu)
            {\chi}( |\eta|\,t )k(\eta,t,s)
            {\chi}( |\eta| \, s )\,\eta^{\mu-1} \, d\eta
            \label{kernel1}\\
  k_f(t,s)&=\sigma_0(t)\,\sigma_1(s)\int_{\calC} f(\eta^\mu)
            {\chi}( |\eta|\,t ) k(\eta,t,s)
            (1-{\chi})( |\eta| \, s )\,\eta^{\mu-1} \, d\eta
            \label{kernel2a}\\
  k_f(t,s)&=\sigma_0(t)\,\sigma_1(s)\int_{\calC} f(\eta^\mu)
            (1-{\chi})( |\eta|\,t)k(\eta,t,s)
            {\chi}( |\eta| \, s )\,\eta^{\mu-1} \, d\eta
            \label{kernel2b}\\
  k_f(t,s)&=\sigma_0(t)\,\sigma_1(s)\int_{\calC} f(\eta^\mu)
            (1-{\chi})( |\eta| \, t)k(\eta,t,s)
            (1-{\chi})( |\eta| \, s )\,\eta^{\mu-1} \, d\eta
            \label{kernel3}
 \end{align}
the associated integral operators satisfies \eqref{normbound}.  
To begin with case \eqref{kernel1} we use the fact 
that, for some fixed $\eps > 0$, 
  $$|{k}(\eta,t,s)|\le c|\eta|^{-\mu+2\eps}
    t^{-\frac{n+1}{2}+\eps}s^{-\frac{n+1}{2}+\eps},$$ 
uniformly in $\eta\in{\calC}$ and $t,s>0$, cf.\ \eqref{h4.B.4} and \eqref{h4.X}. Then
\begin{align*}
  |k_f(t,s)|&\le 
    c\,\| f \|_{\infty}\,\sigma_0(t)\,\sigma_1(s)\,
    t^{-\frac{n+1}{2}+\eps}s^{-\frac{n+1}{2}+\eps}
    \int_1^\infty \varrho^{-1+2\eps}
    {\chi}(\varrho t){\chi}(\varrho s)\,d\varrho\\
  &\le\| f \|_{\infty}\, \frac{c}{2\eps}\,
    \min\left(\frac{1}{t},\frac{1}{s}\right)^{2\eps}\sigma_0(t)\,\sigma_1(s)\,
    t^{-\frac{n+1}{2}+\eps}s^{-\frac{n+1}{2}+\eps}.
\end{align*}
Since the factor $\frac{c}{2\eps}\sigma_0(t)\,\sigma_1(s)$ can be estimated from above 
by a constant, it remains to consider the kernel 
$t^{-\frac{n+1}{2}+\eps}s^{-\frac{n+1}{2}+\eps}\min(\frac{1}{t},\frac{1}{s})^{2\eps}$. 
Because this kernel is symmetric in $s$ and $t$, indeed it suffices to treat
\begin{equation}\label{kern1.4}
  k(t,s)=\begin{cases}
          t^{-\frac{n+1}{2}-\eps}s^{-\frac{n+1}{2}+\eps}&:s\le t\\
          0&:s>t.
         \end{cases}
\end{equation}
Recalling the Hardy inequality 
  $$\int_0^\infty\Big(\int_0^t g(s)\,ds\Big)^pt^{-1-r}\,dt\le 
    \mbox{$\left(\frac{p}{r}\right)^p$}\int_0^\infty g(t)^pt^{p-1-r}\,dt,$$
which holds for any non-negative function $g$ on $\rz_+$ and $r>0$ 
(cf.\ \cite{StWe}, Lemma 3.14, page 196), and denoting by $G$ the integral operator 
with kernel \eqref{kern1.4}, we have
\begin{align*}
  \|Gu\|^p_{\calH^{0,0}_p(X^\wedge)}&\le
     \int_0^\infty\Big(\int_0^\infty 
     k(t,s)|u(s)|\,s^nds\Big)^pt^{\frac{n+1}{2}p-1}\,dt
     =\int_0^\infty\Big(\int_0^t s^{\frac{n-1}{2}+\eps}|u(s)|\,ds\Big)^p t^{-1-p\eps}\,dt\\
  &\le\mbox{$\left(\frac{p}{p\eps}\right)^p$}\int_0^\infty 
     |u(t)|^p\,t^{\frac{n+1}{2}p-1}dt=
     \mbox{$\left(\frac{1}{\eps}\right)^p$}\|u\|^p_{\calH^{0,0}_p(X^\wedge)},
\end{align*}
which completes the proof for the case \eqref{kernel1}. The proofs for the 
cases \eqref{kernel2a}, \eqref{kernel2b} and \eqref{kernel3} can be obtained similarly, 
cf. \cite{CSS1}.
\end{proof}
    
Let us now turn our attention to the analysis of (\ref{mellinterm}). We make use of the 
structure of $\wt{h}$ as described in \eqref{h4.E}, \eqref{h4.F}, \eqref{h4.G}, and 
Definition \ref{h4.3}. Then we have to handle two terms involving 
$\wt{a}_0$ and $\wt{a}_1$, and a term coming from the symbol kernel $\wt{d}$. 

Let us recall that any Mellin symbol 
$h\in M\!S^{\nu}(\rz_+\times\rz^n\times\Gamma_{\beta}\times\rz^n)$, i.e.\ any 
smooth function satisfying, for any order of derivatives,  
 $$|(t\partial_{t})^k\partial^{\alpha^\prime}_{x}\partial^l_{\tau}\partial^\alpha_{\xi}
   h(t,x,\beta+i\tau,\xi)|\le c  \spk{\tau,\xi}^{\nu-l-|\alpha|}, $$
induces continuous operators from 
$\calH^{s,\frac{n+1}{2}-\beta}_p(\rz_+\times\rpbar^n)$ to 
$\calH^{s-\nu,\frac{n+1}{2}-\beta}_p(\rz_+\times\rpbar^n)$.  


\begin{proposition}\label{lower}
Let $\wt{a}_1\in S^{-\mu-1}(\rpbar\times\rz^n\times\rz^{1+n}\times\Sigma)$ 
be compactly supported in $(t,x)$ and
 $$h_f(t,x,\mbox{$\frac{n+1}{2}$}-\gamma+i\tau,\xi)=
   t^\mu\int_{\partial\Sigma} f(\eta^\mu)\,
   \wt{a}_1(t,x,\tau,\xi,t \eta)\,\eta^{\mu-1}\,d\eta$$
with $f\in H$. This defines a symbol 
$h_f\in M\!S^{-1}(\rz_+\times\rz^n\times\Gamma_{\frac{n+1}{2}-\gamma}\times\rz^n)$ 
and the symbol estimates of $\| f \|^{-1}_{\infty}\,h_f$ are uniform in $0\not=f\in H$. 
Consequently, 
  $$\|\text{\rm op}_M^{\gamma-\frac{n}{2}}\op_+(h_f)\|_{\mathcal{L}
     (\calH^{0,\gamma}_p(\rz_+\times\rpbar^n))}\le 
     M_p\, \| f \|_{\infty}\qquad \forall\;f\in H$$ 
for a suitable constant $M_p\ge0$. 
\end{proposition}
\begin{proof}
Without loss of generality, we consider the case $\gamma = \frac{n+1}{2}$. 
We have to show that 
\begin{equation}\label{ms0estimate}
 | \partial^l_{\tau} (t\partial_{t})^k \partial^\alpha_{\xi} \partial^\beta_{x}
 h_f(t,x,i\tau,\xi) | \,
 \langle\tau,\xi \rangle^{1+l+|\alpha|}  \, \| f \|^{-1}_{\infty}
\end{equation}
is uniformly bounded for $t>0$, $x \in \rz^n$, $\tau\in\rz$ and $0\not=f\in H$.
Observing that $t\partial_t t^\mu=\mu t^\mu$ and
  $$t\partial_t\left(\wt{a}_1
    (t,x,i\tau,\xi,t\eta)\right)=
    (t\partial_t\wt{a}_1)(t,x,i\tau,\xi,t\eta)+
    (\eta\partial_\eta\wt{a}_1)(t,x,i\tau,\xi,t\eta),$$
we see that the totally characteristic derivative gives rise to terms
of the same type as $\wt{a}_1$. 
Since the derivatives with respect to $x$, $\tau$ and $\xi$ can be 
taken under the integral sign, we may assume from the beginning that 
$\wt{a}_1\in S^{-\mu-1-j}(\rpbar\times\rz^n\times\rz^{1+n}\times\Sigma)$ 
and show that 
  $$|h_f(t,x,i\tau,\xi)|\le c  \, \langle\tau, \xi \rangle^{-1-j} \, 
    \| f \|_{\infty}$$
uniformly in $t>0$, $x \in \rz^n$ and $f\in H$. By hypothesis, we have
  $$|h_f(t,x,i\tau,\xi)|\le c \, \| f \|_{\infty} \, 
    t^\mu\int_0^{+\infty}\spk{\tau,\xi,t\varrho}^{-\mu-1-j}\,
    \varrho^{\mu-1} \, d\varrho,$$
The transformation $\varrho=t^{-1}\langle \tau,\xi\rangle\sigma$ together with the identity 
$\spk{\tau,\xi,\spk{\tau,\xi}\sigma}=\spk{\tau,\xi}\spk{\sigma}$ yields 
\begin{align*}
 |h_f(t,x,i\tau,\xi)|&\le 
    c \, \| f \|_{\infty} \int_{0}^{+\infty} 
    \spk{\tau,\xi}^{-\mu-1-j}\spk{\sigma}^{-\mu-1-j}\spk{\tau,\xi}^{\mu-1}\sigma^{\mu-1}
    \spk{\tau,\xi}\,d\sigma\\
 &\le c \, \| f \|_{\infty} \spk{\tau,\xi}^{-1-j} 
    \int_{0}^{+\infty}\spk{\sigma}^{-2}\,d\sigma.
\end{align*}
This finishes the proof.
\end{proof}

\begin{proposition}\label{principal}
Let $\wt{a}_0\in S^{-\mu}(\rpbar\times\rz^n\times\rz^{1+n}\times\Sigma)$ be as described 
in Definition {\rm\ref{h4.3}} with compact $(t,x)$-support. Define 
 $$h_f(t,x,\mbox{$\frac{n+1}{2}$}-\gamma+i\tau,\xi)=
   t^\mu\int_{\partial\Sigma} f(\eta^\mu)\,
   \wt{a}_0(t,x,\tau,\xi,t \eta)\,\eta^{\mu-1}\,d\eta$$
with $f\in H$. Then 
$h_f\in MS^{0}(\rz_+\times\rz^n\times\Gamma_{\frac{n+1}{2}-\gamma}\times\rz^n)$, 
and the symbol estimates of $\| f \|^{-1}_{\infty}\,h_f$ are uniform in $0\not=f\in H$. 
Therefore, for a suitable constant $M_p\ge0$, 
  $$\|\text{\rm op}_M^{\gamma-\frac{n}{2}}\op_+(h_f)\|_{\mathcal{L}
     (\calH^{0,\gamma}_p(\rz_+\times\rpbar^n))}\le 
     M_p\, \| f \|_{\infty}\qquad \forall\;f\in H.$$ 
\end{proposition}
\begin{proof}
It is sufficient to show the symbol estimate for $h_f$; the estimates for derivatives 
of $h_f$ are obtained similarly, arguing as in the proof of Proposition \ref{lower}. 
The change of variables $\eta\mapsto t^{-1}\eta$ yields  
\begin{equation}\label{h5.Z}
 h_f(t,x,\mbox{$\frac{n+1}{2}$}-\gamma+i\tau,\xi)=
 \int_{\partial\Sigma} f(t^{-\mu}\eta^\mu)\,
 \wt{a}_0(t,x,\tau,\xi,\eta)\,\eta^{\mu-1}\,d\eta.
\end{equation}
Let us now denote by $\calC_{(\tau,\xi)}$ the natural parametrization of the boundary of 
$\{\eta\in\cz\setminus\Sigma\st |\eta|<c_0\spk{\tau,\xi}\}$, where $c_0$ is associated 
with $\wt{a}_0$ as in Definition \ref{h4.3}. 
As the integrand in \eqref{h5.Z}, for fixed $(t,x,\tau,\xi)$, is holomorphic in $\eta$ 
outside $\Sigma\cup\{|\eta|<c_0\spk{\tau,\xi}\}$ and decays there as 
$|\eta|^{-1-\eps}$ for some $\eps>0$ (recall the decay property of functions in $H$), 
we may replace the integration over $\partial\Sigma$ in \eqref{h5.Z} by integration 
over $\calC_{(\tau,\xi)}$  and obtain  
 $$|h_f(t,x,\mbox{$\frac{n+1}{2}$}-\gamma+i\tau,\xi)|\le c\,\|f\|_\infty\,
   \text{\rm length}(\calC_{(\tau,\xi)})\,\spk{\tau,\xi}^{-\mu}\,\spk{\tau,\xi}^{\mu-1}
   \le c\,\|f\|_\infty,$$
since the length of $\calC_{(\tau,\xi)}$ is less than $(2+2\pi)c_0\spk{\tau,\xi}$. 
\end{proof}


\begin{proposition}\label{singg}
Let ${d}$ be a symbol kernel as in \eqref{h4.F}, \eqref{h4.G}, with  
compact support in $(t,x^\prime)$. Define, for $f\in H$, the symbol kernel  
 $$h_f(t,x^\prime,\mbox{$\frac{n+1}{2}-\gamma+i\tau$},\xi^\prime;u,v) = 
   t^\mu\int_{\partial\Sigma} f(\eta^\mu) \, 
   {d}(t,x^\prime,\tau,\xi^\prime;u,v)\,\eta^{\mu-1}\,d\eta.$$ 
Then there exists a constant $M_p\ge0$ such that 
\begin{equation}\label{h5.W}
 \|\opm{\gamma-\frac{n}{2}}\op^\prime
 (h_f)\|_{\calL(\calH^{0,\gamma}_p(\rz_+\times\rpbar^n))}
 \le M_p\, \| f \|_{\infty}\qquad \forall\;f\in H.
\end{equation}
\end{proposition}
\begin{proof}
Without loss of generality, $\gamma=\frac{n+1}{2}-\frac{1}{p}$. For this choice  
$\calH^{0,\gamma}_p(\rz_+\times\rpbar^n)=L_p(\rz_+\times\rz^n_+)$. 
We are now going to show that 
\begin{equation}\label{h5.Y}
 |\partial^k_\tau \partial^\alpha_{\xi^\prime} (t\partial_t)^l \partial^\beta_{x^\prime}
 h_f(t,x^\prime,\mbox{$\frac{1}{p}$}+i\tau,\xi^\prime;u,v)|
 \le c\,\|f\|_{\infty} \, \frac{\spk{\tau,\xi^\prime}^{-k-|\alpha|}}{u+v}
\end{equation}
for any $k,l\in\nz_0$ and $\alpha,\beta\in\nz_0^{n-1}$. 
It is enough to consider the case $k=l=0$ and $\alpha=\beta=0$, since the terms for 
higher order derivatives are of the same kind, cf.\ the proof of Proposition \ref{lower}. 

Inserting the explicit form of ${d}$, cf.\ \eqref{h4.F}, we obtain 
\begin{equation*}
 h_f(t,x^\prime,\mbox{$\frac{1}{p}$}+i\tau,\xi^\prime;u,v)=
 t \int_{\partial\Sigma} f(\eta^\mu) \, 
 \wt{d}(x^\prime,\xi^\prime,t\eta;[\xi^\prime,\tau,t\eta] u, 
 [\xi^\prime,\tau,t\eta] v)\,(t\eta)^{\mu-1}\,d\eta.
\end{equation*}
Since $[\cdot]\sim\spk{\cdot}$ and $d$ is rapidly decreasing in $(u,v)$, hence in $u+v$, 
we can estimate 
\begin{align*}
  |h_f(t,x^\prime,\mbox{$\frac{1}{p}$}+i\tau,\xi^\prime;u,v)|
& \le
  c \|f\|_{\infty}\int_0^\infty (\varrho t)^{\mu-1} \, 
  \spk{(u+v)\spk{\xi^\prime,\tau,\varrho t}}^{-2} \,
  \spk{\xi^\prime,\tau,\varrho t}^{1-\mu} \, t \, d\varrho\\
& = c \|f\|_{\infty} \int_0^\infty \sigma^{\mu-1} \, 
  \spk{(u+v)\spk{\tau,\xi^\prime}\spk{\sigma}}^{-2} \, \spk{\sigma}^{1-\mu}\,
  \spk{\tau,\xi^\prime}\,d\sigma;
\end{align*}
for the last identity we made use of the change of variables 
$\varrho=t^{-1}\spk{\tau,\xi^\prime}\sigma$. The change of variables 
$r=(u+v)\spk{\tau,\xi^\prime}(1+\sigma)$ then yields 
\begin{align*}
  |h_f(t,x^\prime,\mbox{$\frac{1}{p}$}+i\tau,\xi^\prime;u,v)|
& \le c \|f\|_{\infty} \int_0^\infty \sigma^{\mu-1} \, 
  \spk{(u+v)\spk{\tau,\xi^\prime}(1+\sigma)}^{-2} \, \spk{\sigma}^{1-\mu}\,
  \spk{\tau,\xi^\prime}\,d\sigma\\
& \le c \|f\|_{\infty} \frac{1}{u+v}\int_0^\infty \spk{r}^{-2}\,dr, 
\end{align*}
i.e. \eqref{h5.Y}. Using the continuity of the Hilbert transform, 
 $$u\mapsto Hu=\int_0^\infty\frac{u(s)}{\cdot+s}\,ds:
   L_p(\rz_+)\longrightarrow L_p(\rz_+),$$
and the $L_p(\rz^n)$-continuity of standard zero order pseudodifferential operators, 
assertion \eqref{h5.W} is obtained as follows: Let us write for short $y=(t,x^\prime)$ 
and $\op_{y}=\op_{M,t}^{\gamma-\frac{n}{2}}\op_{x^\prime}$. Then, 
for $\varphi=\varphi(u,y)\in\cicomp(\rz_+\times\rz_+^n)$,  
 $$\big(\op_M^{\gamma-\frac{n}{2}}\op^\prime(h_f)\varphi\big)(u,y)=
   \int_0^\infty\big(\op_{y}(h_f(u,v))\varphi\big)(v,y)\,dv,$$
and therefore 
\begin{align*}
 \|\op^\prime(h_f)\varphi\|_{L_p(\rz_+\times\rz^{n}_+)}^p&\le 
   \int_0^\infty\int\Big(\int_0^\infty|\big(\op_y(h_f(u,v))\varphi\big)(u,y)|dv\Big)^p\, 
   dydu\\
 &\le \int_0^\infty\Big(\int_0^\infty\Big(
   \int|\big(\op_y(h_f(u,v))\varphi\big)(u.y)|^p\,dy\Big)^{1/p}dv\Big)^p du. 
\end{align*}
The second estimate is due to Minkowski's inequality for integrals. Thus  
\begin{align*}
 \|\op^\prime(h_f)\varphi\|_{L_p(\rz_+\times\rz^{n}_+)}^p&\le 
   c\|f\|_\infty^p\int_0^\infty\Big(\int_0^\infty\frac{1}{u+v}
   \|\varphi(v,y)\|_{L_p(\rz_+\times\rz^n_{+,y})}\,dv\Big)^p\,du\\
 &\le c\|f\|_\infty^p \int_0^\infty \|\varphi(u,y)\|_{L_p(\rz_+\times\rz^n_{+,y})}^p\,du
   =c\|f\|_\infty^p\|\varphi\|_{L_p(\rz_+\times\rz^n_+)}^p.
\end{align*}
This finishes the proof. 
\end{proof}

\section{Parameter-ellipticity of the minimal extension}\label{sectionh6}
In Theorem \ref{hinfty} we showed the existence of a bounded $H_\infty$-calculus for 
a closed extension $\ulA_T$, assuming that the resolvent exists in the sector 
$\Lambda$ and has a suitable structure. An obvious problem is now to find conditions on $A$ 
and $T$ which are more easily checked and which ensure all the required assumptions 
on $\ulA_T$. In this section we shall give such conditions for the case  
$\ulA_T=A_{T,\min}$. In fact, these conditions are obtained by combining the concept of 
parameter-ellipticity in Schulze's cone calculus and the observations from 
Section 3.2 of \cite{CSS2}. 

As described in Section \ref{sectionh4.2}, $\dz$-ellipticity of $\spalte{A}{T}$ is 
characterized by the invertibility of the principal symbol $\sigma_\psi^\mu(A)$, the 
rescaled symbol $\wt\sigma_\psi^\mu(A)$, the boundary symbol 
$\sigma_\partial^\mu\spalte{A}{T}$, and the rescaled boundary symbol 
$\wt\sigma_\partial^\mu\spalte{A}{T}$. We shall now pose 
stronger conditions, which we call {\em parameter-ellipticity} with respect to the sector 
$\Lambda$. The first condition is an analog of condition \eqref{h3.Z}: 
\begin{itemize}
 \item[{\bf (E1)}] {\em Both $\sigma_\psi^\mu(A)$ and $\wt\sigma_\psi^\mu(A)$ 
  pointwise do not have spectrum in $\Lambda$ $($for non-zero covariables$)$}.   
\end{itemize}
For $\spalte{A}{T}$ as in \eqref{h4.A.2} we consider the boundary symbol and 
rescaled boundary symbol 
\forget{
For notational convenience we use in the following the same letters $E$ and $F$ 
also for the restrictions of those bundles to subsets of $\dz$, e.g.\ we write 
$E$ instead of $E|_\bz$. Let $\pi_\bz:T^*\intb\to\intb$ and 
$\pi_{\partial\bz}:T^*\partial\bz\to\partial\bz$ denote the canonical projections. 
Then} 
\begin{equation*}
 \sigma_\partial^\mu\begin{pmatrix}A\\ T\end{pmatrix}:
   \calS(\rpbar,E^\prime)\longrightarrow 
   \begin{pmatrix}
    \calS(\rpbar,E^\prime)\\ \oplus\\ F^\prime
   \end{pmatrix},\qquad
 \wt\sigma_\partial^\mu\begin{pmatrix}A\\ T\end{pmatrix}:
   \calS(\rpbar,E^{\prime\prime})\longrightarrow
   \begin{pmatrix}
    \calS(\rpbar,E^{\prime\prime})\\ \oplus\\ F^{\prime\prime}
   \end{pmatrix}, 
\end{equation*}
where $E^\prime,\,E^{\prime\prime}$ and $F^\prime,\,F^{\prime\prime}$ are the 
corresponding pull-backs of $E$ and $F$ to $T^*\intb$ and $T^*\partial\bz$, respectively. 
We require that 
\begin{itemize}
 \item[{\bf (E2)}] {\em Both 
  $\spalte{\lambda}{0}-\sigma_\partial^\mu\spalte{A}{T}$ and 
  $\spalte{\lambda}{0}-\wt\sigma_\partial^\mu\spalte{A}{T}$ are pointwise 
  invertible on $(T^*\intb\times\Lambda)\setminus\{0\}$ and 
  $(T^*\partial\bz\times\Lambda)\setminus\{0\}$, respectively}. 
\end{itemize}

Mainly to ensure the identity \eqref{h4.A.5}, we pose a condition on the 
so-called {\em conormal symbol}  
$\sigma_M^\mu\begin{pmatrix}A\\ T\end{pmatrix}(z)=
\begin{pmatrix}\sigma_M^\mu(A)(z)\\ \sigma_M^\mu(T)(z)\end{pmatrix}$. 
This is a holomorphic family in $z\in\cz$ of boundary value problems on $X$ obtained 
in the following way: Using the representation of $A$ as in \eqref{fuchs}, 
 $$\sigma_M^\mu(A)(z):=\smsum_{j=0}^\mu a_j(0)z^j$$
is a holomorphic family of differential operators on $X$. Similarly, using the notation 
from \eqref{h4.A.6}, one defines 
$\sigma_M^\mu(T)=(\sigma_M^0(T_0),\ldots,\sigma_M^{\mu-1}(T_{\mu-1}))$ by  
 $$\sigma_M^j(T_j)(z)=\gamma_0\circ\sigma_M^j(B_j)(z).$$
\forget{ 
where $\gamma_0$ denotes restriction to the boundary of $X$. 
}

It can be shown that $\sigma_M^\mu\spalte{A}{T}$ is meromorphically
invertible in case $\spalte{A}{T}$ is $\dz$-elliptic (in the sense of 
Section \ref{sectionh4.2}). Then the condition is that 
\begin{itemize}
 \item[{\bf (E3)}] 
  {\em 
  $\sigma_M^\mu\begin{pmatrix}{A}\\{T}\end{pmatrix}(z):\ci(X,E)\longrightarrow
  \begin{matrix}\ci(X,E)\\ \oplus\\ \ci(\partial X,F)\end{matrix}$ is 
  invertible for each $z$ with $\re z=\frac{n+1}{2}-\gamma-\mu$}. 
\end{itemize}

\begin{remark}\label{h6.1}
The invertibility of the conormal symbol from {\rm(E3)} is equivalent to that of   
 $$\sigma_M^\mu\begin{pmatrix}A\\ T\end{pmatrix}(z):H^s_p(X,E)\longrightarrow
   \begin{matrix}
    H^{s-\mu}_p(X,E)\\ \oplus\\ 
    \mathop{\mbox{\Large$\oplus$}}\limits_{j=0}^{\mu-1} 
    B_{pp}^{s-j-1/p}(\partial X,F_j) 
   \end{matrix}$$ 
for $s\ge 0$ and $1<p<\infty$. In fact, according to a result of 
Grubb \cite[Theorem 1.12]{Grub}, invertibility for one choice of $s$ and $p$ implies 
that the inverse also is an element in Boutet de Monvel's calculus and thus infers 
the invertibility for every other choice. Equivalence with the invertibility on spaces 
of smooth functions then follows from the fact that the kernel and the cokernel of an 
elliptic operator consist of smooth functions. 
\end{remark}

For $A$ as in \eqref{fuchs}, we define the so-called {\em model cone operator} $\wh{A}$ 
on $X^\wedge=\rz_+\times X$ as 
 $$\wh{A}=t^{-\mu}\smsum_{j=0}^\mu a_j(0)(-t\partial_t)^j.$$
Similarly, we define $\wh{T}=(\wh{T}_0,\ldots,\wh{T}_{\mu-1})$ by  
$\wh{T}_j=\gamma_0\circ\wh{B}_j$, cf.\ \eqref{h4.A.6}, where now $\gamma_0$ denotes the 
restriction to $\partial X^\wedge=\rz_+\times\partial X$. For the analysis of $\wh{A}$, 
one uses a special scale of Sobolev spaces $\calK^{s,\gamma}_p(X^\wedge)$ 
on $X^\wedge$ with $s\in\rz$ and $1<p<\infty$, namely 
 $$\calK^{s,\gamma}_p(X^\wedge)=\{u\in H^s_{p,{\rm loc}}(X^\wedge)\st 
   \omega u\in\calH^{s,\gamma}_p(\dz)\text{ and } (1-\omega)u\in 
   H^s_{p,{\rm cone}}(X^\wedge)\};$$
here $\omega$ is a cut-off function located near $t=0$, and the subscript `{\rm cone}' 
indicates that we do not consider $X^\wedge$ with its product structure, but as an 
$SG$-manifold, cf.\ Section 4.2 in \cite{ScSc2} for more details. 
If $E^\wedge$ is the pull-back of $E|_X$ to $X^\wedge$ under the canonical projection 
$X^\wedge\to X$, this definition also extends to sections, i.e.\ 
we may define $\calK^{s,\gamma}_p(X^\wedge,E^\wedge)$. 

Analogously to Section \ref{sectionh4.3}, we then consider $\wh{A}$ as an unbounded 
operator 
\begin{equation}\label{h6.C}
 \wh{A}:\calS^\infty(X^\wedge,E^\wedge)_{\wh{T}}\subset 
 \calK^{0,\gamma}_p(X^\wedge,E^\wedge)\longrightarrow 
 \calK^{0,\gamma}_p(X^\wedge,E^\wedge),
\end{equation}
where $\calS^\infty(X^\wedge,E^\wedge)$ are the smooth sections of $E^\wedge$ that 
vanish rapidly for $t\to\infty$ and vanish to infinite order in $t=0$. 
The main ingredients for the analysis of the closed extensions of $\wh{A}$ are: 
\begin{itemize}
 \item[i)] $\wh{T}$ has a right-inverse that belongs to the cone calculus for boundary 
  value problems for the infinite cone (cf.\ for example \cite{KaSc}).
 \item[ii)] For a fixed $0\not=\lambda_0\in\Lambda$, 
  $\spalte{\lambda_0-\wh{A}}{\wh{T}}$ is an elliptic element 
  in the cone calculus and one can construct a parametrix 
  $\begin{pmatrix}\wh{R}&\wh{K}\end{pmatrix}$ inverting it modulo finite rank 
  operators and such that $\wh{T}\wh{R}=0$. 
\end{itemize}
For the analysis of closed extensions of $A$ on $\dz$ both the corresponding right-inverse 
as well as the parametrix were constructed in Lemma 3.4 and Propositions 3.3, 3.7 of 
\cite{CSS2} relying on results of \cite{Grub} for boundary value problems on smooth 
manifolds. Both constructions extend to $X^\wedge$. 

\begin{theorem}\label{h6.2}
Under conditions {\rm(E1)} to {\rm(E3)} the following statements hold: 
\begin{itemize}
 \item[a)] $\calS^\infty(X^\wedge,E^\wedge)_{\wh{T}}$ is a dense subspace of 
  $\calK^{s,\gamma}_p(X^\wedge,E^\wedge)_{\wh{T}}$ for any $1<p<\infty$ and 
  $s,\gamma\in\rz$ with $s>\mu-1+\frac{1}{p}$. 
 \item[b)] The domain of the closure of $\wh{A}$ from \eqref{h6.C}, which we denote by 
  $\wh{A}_{\wh{T},\min}$,  coincides with 
  $\calK^{\mu,\gamma+\mu}_p(X^\wedge,E^\wedge)_{\wh{T}}$. 
\end{itemize}
\end{theorem} 
\begin{proof}
a) Using i) above, we obtain a projection 
$P:\calK^{s,\gamma}_p(X^\wedge,E^\wedge)\to\calK^{s,\gamma}_p(X^\wedge,E^\wedge)_{\wh{T}}$ 
within the cone calculus on the infinite cone $X^\wedge$ 
(cf.\ Section 2.2.3 in \cite{KaSc}). 
Then one argues as in the proof of Corollary 3.10 in \cite{CSS2}. 

b) The continuity of $\wh{A}$ together with a) implies that    
$\calK^{\mu,\gamma+\mu}_p(X^\wedge,E^\wedge)_{\wh{T}}\subset\calD(\wh{A}_{\wh{T},\min})$. 
The reverse inclusion follows as in the proof of Proposition 4.2 in \cite{CSS2} with 
the special parametrix from ii). 
\end{proof}

Our next -- and final -- requirement is that 
\begin{itemize}
 \item[{\bf(E4)}] {\em $\wh{A}_{\wh{T},\min}$ does not have spectrum in 
  $\Lambda\setminus\{0\}$}. 
\end{itemize}
As $\wh{A}$ and $\wh{T}$ are invariant under dilations, the spectrum of 
$\wh{A}_{\wh{T},\min}$ is automatically a conical subset of $\cz$. The next 
proposition shows that the spectrum of $\wh{A}_{\wh{T},\min}$ does not depend on 
the choice of $1<p<\infty$: 

\begin{proposition}\label{h6.3}
Assume that conditions {\rm(E1)} to {\rm(E3)} hold and fix a $\lambda_0\not=0$. 
Let us denote, for the moment, by $A_p$ the minimal extension of $\lambda_0-\widehat{A}$ 
subject to $\wh{T}$ in $\calK^{0,\gamma}_p(X^\wedge,E^\wedge)$. 
Suppose that for some $1<p_0<\infty$ the operator $A_{p_0}$ is invertible. 
Then $A_p$ is invertible for all $1<p<\infty$.
\end{proposition}
\begin{proof}
By Theorem \ref{h6.2}, the domain of $A_p$ is 
$\calK^{\mu,\gamma+\mu}_p(X^\wedge,E^\wedge)_{\wh{T}}$. 
Then the invertibility of $A_p$ is equivalent to the invertibility of 
 $$\widehat\calA:=
   \begin{pmatrix}{\lambda_0-\widehat A}\\{\widehat T}\end{pmatrix}:
   \calK^{\mu,\gamma+\mu}_p(X^\wedge,E^\wedge)\longrightarrow 
   \begin{matrix}
     \calK^{0,\gamma}_p(X^\wedge,E^\wedge)\\ \oplus\\ 
     \mathop{\mbox{\Large$\oplus$}}\limits_{j=0}^{\mu-1}
     \calB_{pp}^{\mu-j-1/p,\gamma+\mu-j-1/2}(\partial X^\wedge,F_j^\wedge)
   \end{matrix},$$
see e.g.\ Corollary 7.2 in \cite{CSS2}. Since $\widehat\calA$ is an 
elliptic element in the cone calculus on $X^\wedge$, 
we find a parametrix $\calB$ to $\widehat\calA$ such that 
\begin{eqnarray}\label{linksrechts}
\widehat\calA\calB=I+\calR_1\quad \text{and}\quad\calB\widehat\calA=I+\calR_2,
\end{eqnarray} 
where $\calR_1$ and $\calR_2$ are operators of order $-\infty$ and types $0$ and $\mu$, 
respectively. They have the following mapping properties:
\begin{equation*}
\calR_1:
\begin{matrix}
 \calK^{s,\gamma}_p(X^\wedge,E^\wedge)\\ \oplus\\ 
 \mathop{\mbox{\Large$\oplus$}}\limits_{j=0}^{\mu-1}
 \calB^{s+\mu-j-1/p,\gamma+\mu-j-1/2}_{pp}(\partial X^\wedge,F_j^\wedge)
\end{matrix}
\longrightarrow 
\begin{matrix}
 \calS^{\gamma}(X^\wedge,E^\wedge)\\ \oplus\\
 \mathop{\mbox{\Large$\oplus$}}\limits_{j=0}^{\mu-1}
 \calS^{\gamma+\mu-j-1/2}(\partial X^\wedge,F_j^\wedge) 
\end{matrix}
\end{equation*}
for any $1<p<\infty$ and $s>-1+\frac{1}{p}$, and
 $$\calR_2:\calK^{s,\gamma+\mu}_p(X^\wedge,E^\wedge)\longrightarrow 
   \calS^{\gamma+\mu}(X^\wedge,E^\wedge),
   \qquad 1<p<\infty,\;s>\mu-1+\mbox{$\frac{1}{p}$}.$$ 
In fact, in Section 2.1.6 of \cite {KaSc} these mapping properties are shown for the 
case $p=2$; for the extension to arbitrary $p$ we use \cite{GrKo}. In case $p=p_0$, 
multiplying \eqref{linksrechts} by $\widehat\calA^{-1}$ yields that 
\begin{equation}\label{inv}
 \widehat\calA^{-1}=\calB+\calB\calR_1-\calR_2\widehat\calA^{-1}\calR_1.
\end{equation} 
For each $1<p<\infty$ the right hand side extends to a bounded map 
 $$\begin{matrix}
   \calK^{0,\gamma}_p(X^\wedge,E^\wedge)\\ \oplus\\ 
   \mathop{\mbox{\Large$\oplus$}}\limits_{j=0}^{\mu-1}
   \calB^{\mu-j-1/p,\gamma+\mu-j-1/2}_{pp}(\partial X^\wedge,F_j^\wedge)
   \end{matrix}
   \longrightarrow \calK^{\mu,\gamma+\mu}_p(X^\wedge,E^\wedge)$$
and moreover restricts to a continuous map 
 $$\begin{matrix}
   \calS^{\gamma}_0(X^\wedge,E^\wedge)\\ \oplus\\ 
   \mathop{\mbox{\Large$\oplus$}}\limits_{j=0}^{\mu-1}
   \calS^{\gamma+\mu-j-1/2}_0(\partial X^\wedge,F_j^\wedge)
   \end{matrix}
   \longrightarrow\calS^{\gamma+\mu}_0(X^\wedge,E^\wedge).$$
By density, the right hand side of \eqref{inv} therefore furnishes an inverse to 
$\widehat\calA$ for arbitrary $p$.
\end{proof}

\begin{theorem}\label{thm:himin}
Let $\spalte{A}{T}$ satisfy the conditions \text{\rm (E1)} to 
\text{\rm (E4)}. Then 
 $$A_{T,\min}:\calH^{\mu,\gamma+\mu}_p(\dz,E)_T\subset
   \calH^{0,\gamma}_p(\dz,E)\longrightarrow\calH^{0,\gamma}_p(\dz,E)$$
fulfills the assumptions of Theorem {\rm\ref{hinfty}}. Hence there exists a 
$c\ge0$ such that $c+A_{T,\min}$ has a bounded $H_\infty$-calculus with respect to 
$\cz\setminus\Lambda$. 
\end{theorem}
\begin{proof}
Let us choose parameter-dependent order reductions in the cone algebra on $\bz$, 
 $$R_j(\eta)\in C^{\mu-j}(\bz,\Sigma;\gamma+\mu-j,\gamma,\theta;F_j,F_j),
   \qquad j=0,\ldots,\mu-1,$$
where the $F_j$ are the bundles from \eqref{h4.A.6} and $\theta\in\nz$ is arbitrary. 
Let $R(\eta)=\text{\rm diag}(R_0(\eta),\ldots,R_{\mu-1}(\eta))$. 
Now the conditions (E1) to (E4) are chosen in such a way that 
 $$\calA(\eta):=\begin{pmatrix}1&0\\0&R(\eta)\end{pmatrix}
   \begin{pmatrix}\eta^\mu-A\\T\end{pmatrix}=
   \begin{pmatrix}\eta^\mu-A\\R(\eta)T\end{pmatrix}\in
   C^{\mu,\mu}(\dz,\Sigma;\gamma+\mu,\gamma,\theta;E;E,F)$$
with $F:=F_0\oplus\ldots\oplus F_{\mu-1}$ is a parameter-elliptic element the 
cone calculus for boundary value problems. It follows that there exists a parametrix 
$\calB(\eta)=\begin{pmatrix}B(\eta)&K(\eta)\end{pmatrix}$ and that, for $p=2$, 
\begin{equation}\label{h6.A}
 \calA(\eta):\calH^{\mu,\gamma+\mu}_p(\dz,E)\longrightarrow 
   \begin{matrix}
    \calH^{0,\gamma}_p(\dz,E)\\ \oplus\\ 
    \mathcal{B}^{-\frac{1}{2},\gamma-\frac{1}{2}}_{pp}(\bz,F)
   \end{matrix}
\end{equation}
is bijective for sufficiently large $|\eta|$, and the inverse coincides with the 
parametrix. Also it follows that $B(\eta)\in{C}^{-\mu,0}(\dz,\Sigma;\gamma;E)$ in 
the sense of Section \ref{sectionh4.4.3}. As we shall show below even 
$B(\eta)\in C^{-\mu,0}_{\rm hol}(\dz,\Sigma;\gamma;E)$. 
Next, we note that invertibility of \eqref{h6.A} is equivalent to the invertibility of 
\begin{equation}\label{h6.B}
 \eta^\mu-A:\calH^{\mu,\gamma+\mu}_p(\dz,E)_T\longrightarrow
 \calH^{0,\gamma}_p(\dz,E);
\end{equation}
the inverse of \eqref{h6.B} is just $B(\eta)$. By Theorem \ref{h4.2.4} we can conclude 
that the invertibility of \eqref{h6.B} then also is true for arbitrary $1<p<\infty$, 
and the inverse again coincides with $B(\eta)$. 

Now let $B(\eta)$ be as in \eqref{h4.D}. Decomposing $\wt{h}$ as in \eqref{h4.D.5} and  
\eqref{h4.E} yields local symbols 
$\wt{a}\in S^{-\mu}_{\cl}(\rpbar\times\rz^n_x\times\rz^{1+n}\times\Sigma)$. 
By parametrix construction, the leading term is given by inversion of the 
parameter-dependent principal symbol, i.e.\ $\wt a=\wt a_0+\wt a_1$ with  
$\wt{a}_1(t,x,\tau,\xi,\eta)\in 
S^{-\mu-1}_\cl(\rpbar\times\rz^n_x\times\rz^{1+n}_{(\tau,\xi)}\times\Sigma)$ and 
 $$\wt{a}_0(t,x,\tau,\xi,\eta)=\chi(|(\tau,\xi,\eta)|)
   (\eta^{\mu}-\wt\sigma_\psi^\mu(A)(t,x,\tau,\xi))^{-1},$$
where $\chi$ is a 0-excision function. By ellipticity assumption (E1), 
$(\eta^{\mu}-\wt\sigma_\psi^\mu(A)(t,x,\tau,\xi))^{-1}$ is defined for 
$0\not=(\tau,\xi,\eta)\in\rz^{1+n}\times\Sigma$. By homogeneity, it is clear that 
 $$\text{\rm spec}(\wt\sigma_\psi^\mu(A)(t,x,\tau,\xi))\subset \{\lambda\in\cz\st 
   \mbox{$\frac{1}{c_0}$}|(\tau,\xi)|^\mu<|\lambda|<c_0|(\tau,\xi)|^\mu\}$$
for a suitable constant $c_0>1$. Thus
$(\eta^{\mu}-\wt\sigma_\psi^\mu(A)(t,x,\tau,\xi))^{-1}$ is defined on 
 $$\Omega:=\{(\tau,\xi,\eta)\st (\tau,\xi,\eta)\in(\rz^{1+n}\times\Sigma)\setminus\{0\}
   \text{ or }|\eta|\ge c_0|\tau,\xi|\}.$$
As $(\eta^{\mu}-\wt\sigma_\psi^\mu(A)(t,x,\tau,\xi))^{-1}$ is positively homogeneous 
of degree $-\mu$ in $(\eta,\tau,\xi)\in\Omega$, and $\Omega\cap\partial U_1(0)$ is 
compact ($U_1(0)$ denoting the unit ball), the estimates \eqref{h4.H} with $\nu=-\mu$ 
follow. 
\end{proof}

\section{Example: The Dirichlet and Neumann problems for the Laplacian}\label{sectionh7}
We equip $\dz$ with a straight conical metric, i.e. a metric that coincides 
 with $dt^2+t^2g$  on  ${]0,1[}\times X$ for a fixed  metric $g$ on $X$. 
The associated scalar Laplacian $-\Delta$ is a Fuchs-type operator. 
Near $t=0$ it can be written in the form
\begin{equation}\label{eq:delta}
- \Delta=-t^{-2}\,\{(t\partial_t)^2+(n-1)t\partial_t+\Delta_X\}, 
    \qquad n=\text{\rm dim}\,X,
\end{equation}
where $\Delta_X$ denotes the Laplacian on $X$ with respect to $g$. We let 
 $$
 \calA_D := \begin{pmatrix}-\Delta \\ \gamma_0\end{pmatrix}
 \mbox{ and }
 \calA_N := \begin{pmatrix}-\Delta \\ t^{-1} \gamma_1\end{pmatrix}
  $$
be the Dirichlet and Neumann boundary value problems for $-\Delta$, respectively. 
We denote by $-\Delta_{D}$ and $-\Delta_N$ the unbounded operators in  
$\calH^{0,\gamma}_p(\dz)$, acting as $-\Delta$ on the domains 
$\cii(\dz)_{\gamma_0}$ and $\cii(\dz)_{\gamma_1}$, respectively. 
In the sequel we will  write $-\Delta_{D/N}$ to address both operators. 
It is easy to see that both satisfy the ellipticity conditions (E1) and (E2).

Given a function space $\calF$ we will use the notation $\calF_D$ and $\calF_N$ 
in place of $\calF_{\gamma_0}$ and $\calF_{\gamma_1}$ to denote the closed subspace 
of $\calF$ where $\gamma_0$ and $\gamma_1$, respectively, vanish.

\forget{
It was shown in \cite{CSS2} that their adjoint problems are 
$$
\calA_D^* = \begin{pmatrix}-\Delta \\ -i\gamma_0\end{pmatrix}
\mbox{ and }
\calA_N^* = \begin{pmatrix}-\Delta \\ -i t^{-1} \gamma_1\end{pmatrix}.
$$
We denote by $-\Delta_{D}$ and $-\Delta_N$ the unbounded operators in  
$\calH^{0,\gamma}_p(\dz)$, acting as $-\Delta$ on the domains 
$\cii(\dz)_{\gamma_0}$ and $\cii(\dz)_{\gamma_1}$, respectively. 
In the sequel we will  write $-\Delta_{D/N}$ to address both operators. 
Given a function space $\calF$ we will also use the notation
$\calF_D$ and $\calF_N$ in place of $\calF_{\gamma_0}$ and $\calF_{\gamma_1}$
to denote the closed subspace of  $\calF$ where $\gamma_0$ 
and $\gamma_1$, respectively, vanish.
}

\subsection{Closed extensions of the Dirichlet and Neumann Laplacians on $\dz$}

According to the definition in Section \ref{sectionh6}, the 
principal conormal symbols of $\calA_D$ and $\calA_N$ are 
  $$\sigma^2_M(\calA_D)(z)=
  \begin{pmatrix}-z^2+(n-1)z-\Delta_X
  \\
  \gamma_0
  \end{pmatrix}
  \mbox{ and }
 \sigma^2_M(\calA_N)(z)=
  \begin{pmatrix}-z^2+(n-1)z-\Delta_X
  \\
  \gamma_1
  \end{pmatrix}.
  $$
They are invertible, unless $-z^2 + (n-1)z =\lambda_j$ for one of the eigenvalues 
$\lambda_0>\lambda_1>\ldots$ of the boundary problems $\Delta_{X, \gamma_0}$ and 
$\Delta_{X, \gamma_1}$, respectively (recall that $\lambda_0 < 0$ for the  Dirichlet 
problem, while $\lambda_0=0$ for the Neumann problem). This is the case for 
 $z=q_j^+$ or $z=q_j^-$ with 
  \begin{equation*}
   q_j^\pm=\mbox{$\frac{n-1}{2}\pm
   \sqrt{\big(\frac{n-1}{2}\big)^2-\lambda_j}$},
   \qquad j\in\nz_0.
  \end{equation*}
We shall now study the minimal extension of $\Delta_{D/N}$ in 
$\calH^{0,\gamma}_p(\dz)$. We shall require that
\begin{equation}\label{c1}
1-\sqrt{\left(\mbox{$\frac{n-1}{2}$}\right)^2-\lambda_0} < \gamma < 
-1+\sqrt{\left(\mbox{$\frac{n-1}{2}$}\right)^2-\lambda_0}.
\end{equation} 
Of course, this only makes sense, if
\begin{equation}\label{c2}
\left(\mbox{$\frac{n-1}{2}$}\right)^2-\lambda_0>1.
\end{equation} 

\begin{theorem}\label{thm:dom}
Assume \eqref{c1} and \eqref{c2}. 
For $1<p<\infty$, the minimal and maximal extensions coincide both 
for the Dirichlet and the Neumann Laplacians 
on $\calH^{0,\gamma}_p(\dz)$, and their domain is $\calH^{2,\gamma+2}_p(\dz)_{D/N}$. 
In case $\gamma=0$ and $p=2$, the $($minimal$)$ extension is self-adjoint.  
\end{theorem}
\begin{proof}
This is an immediate consequence of Proposition 6.1 of \cite{CSS2}, since condition 
\eqref{c1} implies that the conormal symbols of both $\calA_{D}$ and $\calA_{N}$ are 
invertible for all $z$ with $\frac{n+1}{2}-\gamma-2 \le \Re z \le \frac{n+1}{2}-\gamma$. 
The self-adjointness is a consequence of Theorem 4.5 in \cite{CSS2}, since the adjoint 
problems of $\calA_{D}$ and $\calA_{N}$ are $\spalte{-\Delta}{-i\gamma_0}$ and 
$\spalte{-\Delta}{-i t^{-1}\gamma_1}$, respectively.  
\end{proof}

\begin{remark}\label{rem:E4}
\begin{itemize}
 \item[(a)] Assumption \eqref{c1} implies that 
  both $-\Delta_D$ and $-\Delta_N$, considered as unbounded operators on 
  $\calH^{0,\gamma}_p(\dz)$, satisfy the ellipticity condition {\rm(E3)}.
 \item[(b)] In the Dirichlet case,  condition \eqref{c2} 
  is always true for $n\ge 3$; depending on $\lambda_0$ $($i.e.\ on $X$ and $g)$, 
  it might also hold for $n=1$ or $n=2$.
  In the Neumann case, condition \eqref{c2} only holds for   $n > 3$.
\end{itemize}
\end{remark}

\subsection{Domains of the model cone operator}
Recall that $\widehat{\Delta}_{D/N,\min}$ denotes the closure in 
$\calK^{0,\gamma}_p(X^\wedge)$ of the model cone operator $\wh{\Delta}_{D/N}$ 
considered with domain $\calS^\infty(X^\wedge)_{D/N}$, while 
 $$\calD(-\widehat{\Delta}_{D/N,\max})=\{u\in\calK^{2,\gamma}_p(X^\wedge)_{D/N} \st 
   -\widehat{\Delta}u \in \calK^{0,\gamma}_p(X^\wedge) \}$$ 
defines the closed operator $\widehat{\Delta}_{D/N,\min}$. Note that always 
$\calK^{2,\gamma+2}_p(X^\wedge)_{D/N}\subset\calD(-\widehat{\Delta}_{D/N,\min})$. 
\forget{
In fact, $\calK^{2,\gamma+2}_p(X^\wedge)_{D/N}$ is invariant under 
multiplication by cut-off functions, since the boundary condition does 
not contain $t$-derivatives. Hence, the statement follows from 
\cite{CSS2}, Proposition 4.2, for functions supported close to the singularity. 
For functions supported away from $\{t=0\}$ we recall that 
$\calK^{2,\gamma+2}_p(X^\wedge)$ coincides with the cone Sobolev 
space $H^2_{p,\textrm{cone}}(X^\wedge)$, 
and we deduce the result as in the standard case, cf.\ Section 1.6 in \cite{Grub}. 
}

\begin{lemma}\label{lemma:model}
Under the assumptions of Theorem {\rm \ref{thm:dom}}, the minimal and 
maximal extensions coincide both for the model Dirichlet and Neumann 
Laplacians on $\calK^{0,\gamma}_p(X^\wedge)$ and
$$
\calD(-\widehat{\Delta}_{D/N,\min})
=\calD(-\widehat{\Delta}_{D/N,\max})=
\calK^{2,\gamma+2}_p(X^\wedge)_{D/N}.
$$
\end{lemma}

\begin{proof}
For $u\in \calD(-\widehat{\Delta}_{D,\max})$ we have 
$u\in\calK^{2,\gamma}_p(X^\wedge)_D$ and 
$-\widehat{\Delta}u\in\calK^{0,\gamma}_p(X^\wedge)$. Then
$(1-\omega)u \in \calK^{2,\infty}_p(X^\wedge)_D$, and 
$-\widehat{\Delta} (1-\omega)u \in \calK^{0,\infty}_p(X^\wedge)$. 
Hence, $(1-\omega)u \in \calD(-\widehat{\Delta}_{D,\max})$, thus
$\omega u \in \calD(-\widehat{\Delta}_{D,\max})$. As we may consider 
$\omega u$ as an element of $\calH^{2,\gamma}_p(\dz)$, and since $\Delta$ and 
$\widehat{\Delta}$ have the same form close to $t=0$, we also have 
$\omega u\in\calD(-\Delta_{D,\max})$. By Theorem \ref{thm:dom}, 
$\omega u\in \calD(-\Delta_{D,\min})=\calH^{2,\gamma+2}_p(\dz)_D$. Therefore, 
$u=\omega u + (1-\omega)u \in \calK^{2,\gamma+2}_p(X^\wedge)_D\subset 
\calD(-\widehat{\Delta}_{D,\min})$. 
We conclude that $\calD(-\widehat{\Delta}_{D,\min})= 
\calD(-\widehat{\Delta}_{D,\max})=\calK^{2,\gamma+2}_p(X^\wedge)_D$, and  
the proof of the statement for $-\widehat{\Delta}_D$ is complete. The argument for
$-\widehat{\Delta}_N$ is the same.
\end{proof}

\begin{corollary}\label{selfadjoint}
For $\gamma=0$ and $p=2$ the closure $\Delta_{D/N,\min}$ is self-adjoint. 
\end{corollary}
\begin{proof}
Let $\wh{\Delta}_{F}$ denote the Friedrichs extension of $-\wh{\Delta}_{D/N}$. 
By construction, one has 
 $$\calD(\widehat{\Delta}_{D/N,\min})\subset 
   \calD(\wh{\Delta}_{F})\subset\calD((\widehat{\Delta}_{D/N,\min})^*).$$ 
Let $u\in\calD((\widehat{\Delta}_{D/N,\min})^*)$ be given. 
Using the special parametrix $\wh{R}$ from ii) in Section \ref{sectionh6} 
(for $\wh{A}=-\wh{\Delta}$, $T=D/N$, and $\lambda_0=-1$), and the fact that 
$R^*(\widehat{\Delta}_{D/N,\min})^*\subset(\wh{\Delta}_{D/N}R)^*$, 
we deduce that $u\in\calK^{2,0}_2(X^\wedge)$. Then, arguing as in the proof of 
Theorem 4.5 in \cite{CSS2} to verify there the identity (4.5), we conclude that 
$u\in\calD(\widehat{\Delta}_{D/N,\max})$. It follows from Lemma \ref{lemma:model} that 
$\wh{\Delta}_{F}=\widehat{\Delta}_{D/N,\min}$. 
\end{proof}

\begin{theorem}\label{thm:hinfty}
Assume \eqref{c1} and \eqref{c2}.
Then, for $1<p<\infty$, both the Dirichlet and Neumann Laplacian fulfill the 
ellipticity conditions {\rm(E1)} to {\rm(E4)} of Section {\rm\ref{sectionh6}}.
\end{theorem} 
\begin{proof}
It was noted above that (E1), (E2), and (E3) hold for 
$-\Delta_{D/N}$. It remains to check (E4). 
By spectral invariance, cf.\ Proposition \ref{h6.3},  we 
may assume $p=2$. 

For $\gamma=0$, the Laplacians $-\widehat{\Delta}_{D/N,\min}$ are self-adjoint 
with domain $\calK^{2,2}_2(X^\wedge)_{D/N}$. 
\forget{
and positive. Hence, the associated Friedrichs extensions are 
selfadjoint. As the Friedrichs extension is between the minimal and 
the maximal, and these coincide by Lemma \ref{lemma:model}, its 
domain is $\calK^{2,2}_2(X^\wedge)_{D/N}$. 
}
For $\lambda\notin\rpbar$, we therefore conclude that
\begin{equation}\label{eq:selfadj}
-\widehat{\Delta}_{D/N}-\lambda \;\colon\; 
\calK^{2,2}_2(X^\wedge)_{D/N} \longrightarrow
\calK^{0,0}_2(X^\wedge)
\end{equation}
is invertible. Clearly, \eqref{eq:selfadj} shows the injectivity of 
$-\widehat{\Delta}_{D/N}-\lambda$ on 
$\calK^{2,\gamma+2}_2(X^\wedge)_{D/N}\subset\calK^{2,2}_2(X^\wedge)_{D/N}$ 
for $\gamma>0$. 

Next suppose $\gamma<0$, $u\in\calK^{2,\gamma+2}_2(X^\wedge)_{D/N}$ and
$u\in\ker \{ -{\widehat\Delta}_{D/N} -\lambda \colon 
\calK^{2,\gamma+2}_2(X^\wedge)_{D/N}\rightarrow \calK^{0,\gamma}_2(X^\wedge)\}$. 
Then $u$ is in the maximal domain of $-{\widehat\Delta}_{D/N} -\lambda$, considered 
as an unbounded operator on $\calK^{0,\gamma+2}_2(X^\wedge)$. 
Since condition (E4) holds, $u$ also is in the associated minimal domain, which is 
$\calK^{2,4+\gamma}_2(X^\wedge)_{D/N}$. We can iterate this argument $j$ times, 
until $\gamma+2j > 0$. Then, we conclude that $u=0$, by the previous step.
Thus $-{\widehat\Delta}_{D/N} -\lambda$ is injective on 
$\calK^{2,\gamma+2}_2(X^\wedge)_{D/N}$ for all $\gamma$ satisfying the hypotheses 
and $\lambda\notin\rpbar$. Finally, we note that the adjoint of
$$
-{\widehat\Delta}_{D/N} -\lambda \; \colon\; 
\calK^{2,\gamma+2}_2(X^\wedge)_{D/N} \longrightarrow
\calK^{0,\gamma}_2(X^\wedge)
$$
with respect to the scalar product of $\calK^{0,0}_2(X^\wedge)$ is
$$
-{\widehat\Delta}_{D/N} -\overline{\lambda} \; \colon\; 
\calK^{2,-\gamma+2}_2(X^\wedge)_{D/N} \longrightarrow
\calK^{0,-\gamma}_2(X^\wedge).
$$
Hence, also the adjoint is injective for $\lambda\notin\rpbar$, so 
$-{\widehat\Delta}_{D/N} -\lambda$
is bijective, as claimed.
\end{proof}

\subsection{Maximal $L_p$ regularity of the Cauchy problem for 
Dirichlet and Neumann Laplacians}

As a consequence of Theorem \ref{thm:hinfty} we get the following 
result on the 
solvability of the Cauchy Problem for the Dirichlet and Neumann 
Laplacians.  

\begin{theorem}\label{cauchy}
Let $\Delta$ be the Laplacian as described above, $1<p<\infty$, and assume \eqref{c1} 
and \eqref{c2}. Then the initial boundary value problems 
\begin{equation}\label{cauchy1}
   u'(\tau)-\Delta u(\tau)=f(\tau),\quad 0\le \tau\le T;\qquad u(0)=0,\quad 
   \gamma_j u = 0,
\end{equation}
$j=0,1$, have  a unique solution 
  $$u\in W^1_r\left([0,T],\calH^{0,\gamma}_p(\dz)\right)\,\cap\,
    L_r\left([0,T],\calH^{2,\gamma+2}_q(\dz)_{D/N}\right)$$
for each
  $$f\in L_r\left([0,T],\calH^{0,\gamma}_p(\dz)_{D/N}\right),\qquad 
    1<r<\infty.$$
Furthermore, $u$, $u'$, and $\Delta u$ depend continuously on $f$. 
\end{theorem} 

\begin{proof}
Solving \eqref{cauchy1} is equivalent to solving 
$v'(\tau)-(\Delta_{D/N}-c)v(\tau)=e^{c\tau}f(\tau)$, $v(0)=0$, $\gamma_j v=0$, 
for some $c>0$. 

The operators $-\Delta_{D/N}$ are closed with minimal (and maximal) domain equal to 
$\calH^{2,\gamma+2}_p(\dz)_{D/N}$. According to Theorem \ref{thm:hinfty}, they 
satisfy conditions (E1)-(E4) for each sector $\Lambda$ not containing $\rz_+$. 
Applying Theorem \ref{thm:himin}, we deduce that $-\Delta_{D/N} + c$ has a 
bounded $H_\infty$-calculus for sufficiently large $c$, and  Theorem \ref{cauchy} 
immediately follows from Dore and Venni's theorem, cf.\ Theorem 3.2 in \cite{DoVe}.
\end{proof}

Nazarov \cite{Naza} has studied the Dirichlet and the Neumann problem for the Laplacian 
on infinite cones and wedges in Euclidean space. He shows results on maximal regularity 
in weighted Sobolev spaces using the explicit Green's function. He obtains restrictions 
on the weight which are similar to \eqref{c1}. They do not coincide, however, since he 
works on a different scale of spaces. 

\begin{small}
\bibliographystyle{amsalpha}

\end{small}
\end{document}